\documentclass[11pt,letterpaper]{amsart}
\usepackage[utf8]{inputenc}
\usepackage{amsmath,amssymb,graphicx,setspace,verbatim,url}
\usepackage{color}
\usepackage[colorlinks=true, pdfstartview=FitV, linkcolor=blue, citecolor=blue, urlcolor=blue]{hyperref}
\usepackage{multicol}
\usepackage{enumitem}
\usepackage{float}
\usepackage{soul}
\usepackage{tikz}
\usepackage[normalem]{ulem}
\usepackage[margin=3cm]{geometry}

\usepackage{tikz-cd}
\usetikzlibrary{arrows.meta}
\newenvironment{red}{\relax\color{red}}{\hspace*{.5ex}\relax}
\newcommand{\bcp}{\begin{red} }
\newcommand{\ecp}{\end{red}}
\newcommand{\llN}{\langle\langle}
\newcommand{\rrN}{\rangle\rangle}

\definecolor{darkred}{rgb}{0.7,0,0} % darkred color
\newcommand{\defn}[1]{{\color{darkred}\emph{#1}}} % emphasis of a definition

\usepackage{comment}

\tikzcdset{
  arrow style={tikz,diagrams={>=Triangle}}
}

\input xy
\xyoption{all}

\newtheorem{prop}{Proposition}[section]
\newtheorem{coro}[prop]{Corollary}
\newtheorem{thm}[prop]{Theorem}
\newtheorem{lemma}[prop]{Lemma}
\newtheorem{conjecture}[prop]{Conjecture}

\theoremstyle{definition}
\newtheorem{definition}[prop]{Definition}

\theoremstyle{remark}
\newtheorem{example}[prop]{Example}
\newtheorem{remark}[prop]{Remark}

\DeclareMathOperator{\rank}{rank}

\newcommand{\abs}[1]{\lvert #1 \rvert}
\newcommand{\iso}{\cong}

\newcommand{\GG}{\mathcal{G}}

\newcommand{\HH}{\mathcal{H}}
\newcommand{\KK}{\mathcal{K}}
\newcommand{\LL}{\mathcal{L}}

\newcommand{\Sh}{\mathrm{Sh}}
\newcommand{\Aut}{\mathrm{Aut}}

\newcommand{\Cyc}{\mathrm{C}}

\newcommand{\C}{\mathcal{C}}
\newcommand{\CC}{\mathbb{C}}
\newcommand{\RR}{\mathbb{R}}
\newcommand{\ZZ}{\mathbb{Z}}

\newcommand{\ii}{\boldsymbol{i}}

\usepackage[colorinlistoftodos]{todonotes}

\begin{document}
\title[Coset Geometries and Tameness of Certain Classes of Shephard Groups]{Coset Geometries and Tameness of Certain Classes of Shephard Groups}

\author[Claudio Alexandre Piedade]{Claudio Alexandre Piedade}
\address[Claudio Alexandre Piedade]{Department of Mathematics, Université Libre de Bruxelles, Brussels, Belgium, Orcid number 0000-0002-0746-5893
}
\email{claudio.piedade@ulb.be}

\author[T. Scrimshaw]{Travis Scrimshaw}
\address[T. Scrimshaw]{Department of Mathematics, Hokkaido University, Sapporo  060-0808, Japan}
\email{tcscrims@gmail.com}

\author[Philippe Tranchida]{Philippe Tranchida}
\address[Philippe Tranchida]{Department of Mathematics, Université Libre de Bruxelles, Brussels, Belgium
}
\email{tranchida.philippe@gmail.com}

\setlength{\marginparwidth}{2.5cm}

{}%{}

\thanks{\textit{Declaration of interests:} none}
\date{\today}

\begin{abstract}
%Shephard groups generalise both Coxeter and Artin--Tits groups.
%While the latter have been extensively studied, with numerous results about their algebraic structure and associated geometric objects, much less is known about Shephard groups, and most existing work focuses on finite Shephard groups.
%In this paper, we investigate the structure of the standard parabolic subgroups of Shephard group and of their intersections, especially for low rank cases. In particular, we show that whenever all the standard parabolic subgroups are also Shephard groups, their intersection is also a standard parabolic subgroup. 
%Additionally, we study properties of the natural coset geometries associated with Shephard groups, particularly for finite Shephard groups and right-angled Shephard groups.
%Finally, we introduce a new representation of Shephard groups, mixing the Bureau and the geometric representations, that is applicable to all Shephard groups, even when some of the defining generators have infinite order.
For the class of Shepard groups (which contain Coxeter groups, Artin--Tits groups, and many finite complex reflection groups), we investigate the structure of the coset incidence systems defined using the standard parabolic subgroups.
For the class of tame Shephard groups, which is when all the standard parabolic subgroups are equal to the corresponding Shephard groups, then the intersection of standard parabolic subgroups is also a standard parabolic subgroup.
We show that when the Shephard group is either finite, right-angled, or rank two, then it is tame; this is further evidence that all Shephard groups are tame.
Finally, we mix the Burau and the geometric representations to produce a representation for all Shephard groups, but we show it is not faithful in rank two unless the Shephard group is finite or an Artin--Tits group.
\end{abstract}

\maketitle

\section{Introduction}

Let $\Gamma$ be a simple graph with finite vertex set $I$ and edges $\{i,j\}$ labeled with integers $r_{i,j}=r_{j,i}\in \mathbb{N}_{\geq 3}\cup \{\infty\}$.
The Coxeter group associated with $\Gamma$ is the group $W := W_\Gamma$ is defined by the group presentation
\[
W_\Gamma = 
\left\langle
 s_i, i\in V(\Gamma)
\;\middle\vert\;
\begin{matrix}
s_i^2= 1\textnormal{ for all }i\in I,\\
(s_is_j)^{r_{i,j}}=1\textnormal{ for }\{i,j\}\in E(\Gamma)\textnormal{ and }r_{i,j}\neq \infty\\
(s_is_j)^2=1\textnormal{ if }\{i,j\}\notin E(\Gamma)
\end{matrix}
\right\rangle.
\]
They are named in honor of H.M.S.~Coxeter. He studied these groups from the point of view of reflection groups and polytopes.
Later, Jacques Tits further developed the general theory of Coxeter groups in~\cite{Tits1961}.
However, these groups go beyond polytopes, and have implications in different areas of mathematics, such as Tits buildings, Lie theory, and Springer theory.
The literature on Coxeter groups is quite extensive, and we refer the reader to~\cite{Bourbaki2006,Davis2025,Humphreys1990} for more information.

Coxeter groups are known to be well-behaved in many regards. For example, they have solvable word and conjugacy class problems thanks to their association with root systems (in particular, every Coxeter group can be realized as a reflection group).
In this paper, we will focus on the following properties of Coxeter groups.
For any subset $J\subseteq I$, the group $W_J := \langle s_i \mid i\in J\rangle$ is called a standard parabolic subgroup of $W$.
Here are the three crucial properties of Coxeter groups we will consider (see, e.g.,~\cite[Chapter 4, §1]{Bourbaki2006}):
\begin{enumerate}
    \item $W_J \iso W_{\Gamma_J}$ for every $J\subseteq I$, where $\Gamma_J \subseteq \Gamma$ is the induced subgraph with vertex set $J$;
    \item $W_J \cap W_K = W_{J\cap K}$ for every $J,K\subseteq I$;
    \item $W_J W_K \cap W_J W_L = W_J (W_K \cap W_L)$ for every $J,K,L\subseteq I$.
\end{enumerate}

There is another well-studied and important family of groups, also introduced by Tits~\cite{Tits1966}. These groups are intimately related with Coxeter groups and are called Artin--Tits groups (or sometimes simply Artin groups). They can be seen as a generalization of braid groups.
Just as for a Coxeter group, an Artin--Tits group $A_\Gamma$ is defined from the graph $\Gamma$ by the group presentation
\[
A_\Gamma := 
\left\langle
\begin{matrix}
    \\
 s_i, i\in V(\Gamma)   \\
    \\
\end{matrix}
\right\rvert
\left. \begin{matrix}
\Delta_{s_i,s_j}^{r_{i,j}}=\Delta_{s_j,s_i}^{r_{i,j}}\textnormal{ for }\{i,j\}\in E(\Gamma)\textnormal{ and }r_{i,j}\neq \infty\\
s_is_j=s_js_i\textnormal{ if }\{i,j\}\notin E(\Gamma)
\end{matrix}
\right\rangle.
\]
Notice that, in the definition of a Coxeter group, since $s_i^2 = 1$, each Coxeter relation $(s_is_j)^{r_{i,j}} = 1$ can be written as $\Delta_{s_i,s_j}^{r_i,j} = \Delta_{s_j,s_i}^{r_i,j}$, where
\begin{equation*}
    \Delta_{s_i,s_j}^{r_{i,j}} =
    \begin{cases}
        \underbrace{s_is_js_i\cdots s_j }_{r_{i,j}\textnormal{ factors}} & \textnormal{ if } r_{i,j} \textnormal{ is even,} \\
        \underbrace{s_is_js_i\cdots s_i }_{r_{i,j}\textnormal{ factors}} & \textnormal{ if } r_{i,j} \textnormal{ is odd.}
    \end{cases}
\end{equation*}
With that in mind, the only difference in the presentation of a Coxeter group and an Artin--Tits group are the relations $s_i^2 = 1$, and hence there is a natural projection from $A_{\Gamma} \to W_{\Gamma}$ sending generators to generators.
A striking fact about Artin--Tits groups is that they also satisfy the same three properties $(a)-(c)$ of Coxeter groups described above, where the subgroups $W_J$ are replaced by the subgroups $A_J := \langle s_i\mid i\in J\rangle$ (for more details, see~\cite{van1983homotopy,ParisGodelle2012}).
There are several open problems on Artin--Tits groups, particularly the word problem, the conjugacy problem, and the $K(\pi,1)$ conjecture~\cite{Godelle2012,paris2014k}, which have generated significant interest and research from multiple perspectives.
%\claudio{Travis, is there something we should say here about the Burau representation? To motivate a bit the section on the representation of Shephard groups?}

While studying many of these groups, J.~Tits introduced a general construction of geometrical structure built directly from a group~\cite{Tits1957}.
This structure, nowadays called coset incidence system, takes a group $G$ and a collection of its subgroups $(G_i)_{i\in I}$ as input and outputs an incidence system~\cite{buekenhout2013diagram}.
Geometrical properties of this structure can then be related to algebraic properties of the group $G$, like the above-mentioned properties of Coxeter and Artin--Tits groups, for example.
Indeed, if we consider the Coxeter group $G=W_\Gamma$ and the family of subgroups $(W_{I\setminus\{i\}})_{i\in I}$, the properties $(a)-(c)$ imply that $(W_\Gamma,(W_{I\setminus\{i\}})_{i\in I})$ is a thin residually-connected geometry on which $W_\Gamma$ acts flag-transitively.
If instead we consider the Artin--Tits group $G = A_\Gamma$, we can recover the Artin--Tits flag complex as defined by Paris and Godelle~\cite{ParisGodelle2012}.
Hence, properties $(a)-(c)$ not only provide information about the group, but also about geometrical objects on which these groups act on.

The groups of interest in this paper are Shephard groups.
These appear in the literature with very close connection with Coxeter groups.
They were first described in~\cite{Shephard1952} by Shephard when looking at the automorphism group of regular complex polytopes in $GL_n(\CC)$.
The groups described by Shephard have similar presentations to the ones of Coxeter groups, however generators need not be involutions.
Hence, these groups appear as an interpolation between Coxeter and Artin--Tits groups.
From another point of view, Shephard groups are quotients of Artin--Tits groups, where additional constraints on the order of the generators are imposed.

By the classification of finite complex reflection groups, Shephard and Todd~\cite{ShephardTodd1954} described which finite Shephard groups are subgroups of $GL_n(\CC)$ and, more recently, Goldman~\cite{Goldman2024} proved that all finite Shephard groups are indeed finite complex reflection groups.
It can be easily checked computationally that properties $(a)-(c)$ are true for finite Shephard groups.
Indeed, by considering the natural coset incidence system for this groups, one recovers the incidence structure of the regular complex polytope on which these groups act.

However, when these groups are not finite, the three properties above are not known to hold in general. Indeed, there are very few results on infinite Shephard groups~\cite{2ShepKatherine,Goldman2024,Kapovich1998,piedade2026group}.
Nonetheless, some things are known. Dyer groups, a subclass of Shephard groups, generalise both Coxeter groups and graph products of cyclic groups.
These groups have property all three properties $(a)-(c)$.
Indeed, $(a)$ follows from the seminal paper of Dyer~\cite{Dyer1990}, $(b)$ can be deduced by either~\cite[Corollary 2.6]{Paris2023} or results of this paper (Theorem~\ref{thm:Tame_imply_IP}), while $(c)$ follows from~\cite{piedade2026group,wreathproductPiedadeTranchida} from an iterative construction of these groups by either direct products or free-products with amalgamation.
Note that Dyer groups include right-angled Shephard groups, as these are graph products of cyclic groups.

The goal of this paper is to address these three properties on Shephard groups, particularly properties $(a)$ and $(b)$ which we call, respectively, tameness and intersection property.
In particular, in Section 3 we start by studying Shephard groups of rank $2$. We show these groups satisfy all three properties and describe the gonality of the incidence graph obtained from its coset geometry. Moreover, we give a complete list of all the possible geometries for finite rank two Shephard groups.
In Section 4, we investigate tameness for higher rank Shephard groups. We prove that tame Shephard groups always satisfy the intersection property. Moreover, we improve known results about tameness of rank three Shephard groups.
In Section 5, we give a proof of the three properties for right-angled Shephard groups, independent of the results on Dyer groups. Indeed, we show that any graph product of groups satisfy the three properties, with particular focus on property $(c)$ which gives the flag-transitivity of its natural coset incidence system.
Finally, in Section 6, we give a new representation for Shephard groups that combines the Burau representation for Artin--Tits groups with the geometric representation for Shephard groups. This allows for a geometric representation for Shephard groups with generators of arbitrary orders, even infinite.

\section{Background}\label{sec:prelims}

\subsection{Incidence Geometry and Coset Geometry}
\label{sec:back:subsec:incidence}
    An \defn{incidence system} is a 4-tuple $\Sigma = (X,I,\ast,t)$ consisting of
    \begin{enumerate}
        \item a set $X$ of \defn{elements} of $\Sigma$,
        \item a set $I$ of \defn{types} of $\Sigma$, 
        \item a symmetric and reflexive relation $\ast$ on $X$ called the \defn{incidence relation}, and
        \item a map $t \colon X \to I$, called the \defn{type map} of $\Sigma$, such that distinct elements $x,y \in X$ with $x * y$ satisfy $t(x) \neq t(y)$. 
    \end{enumerate}

%Elements of $t^{-1}(i)$ are called the elements of type $i$ (or $i$-elements) and the set of elements of type $i$ is denoted by $S_i(\Gamma) = S_i$.
The \defn{rank} of $\Sigma$ is the cardinality of the type set $I$.
In an incidence system $\Sigma$, a \defn{flag} is a set of pairwise incident elements.
The rank of a flag $F$ is $\rank(F) = \abs{F}$, the cardinality of $F$.
The type of a flag $F$ is $t(F) = \{t(x) \mid x \in F\}$, the set of types of the elements of $F$. 
%The \defn{cotype} of a flag $F$ is the complement in $I$ of the type of $F$.
A \defn{chamber} is a flag of type $I$.
An incidence system $\Sigma$ is an \defn{incidence geometry} if all its maximal flags are chambers.

The \defn{incidence graph} of $\Sigma$ is a graph with vertex set $X$ with an edge between $x, y \in X$ if and only if $x * y$.
In terms of the incidence graph, the type map is a valid coloring of the graph by $I$, a flag is a clique (isomorphic to a complete graph), and a chamber is a clique with every color.

Let $F$ be a flag of $\Sigma$. An element $x\in X$ is \defn{incident} to $F$ if $x*y$ for all $y\in F$.
The \defn{residue} of $\Sigma$ with respect to $F$, denoted by $\Sigma_F$, is the incidence system formed by all the elements of $\Sigma$ incident to $F$ but not in $F$.
Whenever $\Sigma$ is a geometry, the \defn{rank of a residue} $\Sigma_F$ is defined to be rank$(\Sigma) - \abs{F}$.
If $\Sigma$ is of rank $n$, then the \defn{corank} of a flag $F$ is the difference between the rank of $\Sigma$ and the rank of $F$.

An incidence geometry $\Sigma$ is \defn{connected} if its incidence graph is connected.
It is \defn{residually connected} if all its residues of rank at least two are connected. 
An incidence geometry is \defn{firm}, \defn{thin} or \defn{thick} if all its residues of flags of corank one contain, respectively, at least two elements, exactly two elements, more than two elements.

Let $\Sigma = (X,I,*,t)$ be an incidence geometry.

A \defn{correlation} of $\Sigma$ is a bijection $\phi \colon X \to X$ respecting the incidence relation $*$ and such that, for every $x,y \in X$, if $t(x) = t(y)$ then $t(\phi(x)) = t(\phi(y))$.
If, moreover, $\phi$ fixes the types of every element (i.e.\ $t(\phi(x)) = t(x)$ for all $x \in X$), then $\phi$ is said to be an \defn{automorphism} of $\Sigma$.
The group of all automorphisms of $\Sigma$ is denoted by $\Aut(\Sigma)$. 

If for all $J \subseteq I$ the group $\Aut(\Sigma)$ action is transitive on the set of flags of $\Sigma$ of type $J$, then we say that $\Aut(\Sigma)$ is \defn{flag-transitive} on $\Sigma$.
We say $\Aut(\Sigma)$ is \defn{chamber-transitive} on $\Sigma$ if $\Aut(\Sigma)$ is transitive on the set of chambers of $\Sigma$.
Note that, if $\Sigma$ is a geometry, being chamber-transitive is equivalent to being flag-transitive~\cite[Proposition 2.2]{hypertopes}.

The main interest of this paper is incidence geometries that are obtained from a group $G$ together with a set $(G_i)_{i \in I}$ of subgroups of $G$ as described in~\cite{Tits1957}. 
%\begin{definition}
    A \defn{coset incidence system} $\Sigma(G,(G_i)_{i \in I})$, denoted also as $(G,(G_i)_{i \in I})$, is the incidence system over the type set $I$ as follows.
    \begin{enumerate}
        \item The elements of type $i \in I$ are left cosets $\{g G_i \mid g \in G\}$.
        \item The incidence relation is given by non-empty intersection. More precisely, the element $g G_i$ is incident to the element $k G_j $ if and only if $g G_i \cap k G_j \neq \emptyset$.
    \end{enumerate}
We usually call the subgroups $G_i$ the \defn{(standard) maximal parabolic subgroups}. For $J\subseteq I$, we define \defn{(standard) parabolic subgroups} as $G_J=\bigcap_{i\in J} G_i$, and define $G_{\emptyset} = G$.
When $J=I\setminus\{i\}$, we say that $G_J=G_{I\setminus\{i\}}$ are the \defn{minimal parabolic subgroups}, and we will denote them as $G^i$. Finally, the subgroup $G_I = \bigcap_{i\in I} G_i$ is called the \defn{Borel subgroup} of $\Sigma(G, (G_i)_{i \in I})$. 
    
A coset incidence system that is an incidence geometry is called a \defn{coset geometry}.
Flag-transitive geometries can always be constructed as coset geometries.
The following result adapted from~\cite[Proposition 2.2]{piedade2026group} gives a group theoretical condition to check the flag-transitivity of $G$ on the coset incidence system $\Sigma(G,(G_i)_{i\in I})$. 

\begin{prop}%[{\cite[Proposition 2.2]{piedade2026group}}]
\label{prop:cosetFT}
Let $\Sigma := \Sigma(G,(G_i)_{i\in I})$ be a coset incidence system. Then $G$ is flag-transitive on $\Sigma$ if and only if for every $J,H,K \subseteq I$ we have $G_H G_J \cap G_K G_J = (G_H\cap G_K ) G_J$.
Additionally, if $G$ is flag-transitive on $\Sigma$, then $\Sigma$ is a geometry.
\end{prop}

Notice that $G$ always naturally acts on $\Sigma(G,(G_i)_{i \in I})$ by left multiplication.

Next, we want to consider groups $G$ with a distinguished generating set $S$ and construct a special coset incidence system for the pair $(G, S)$.
Let us discuss general structural properties with regard to subgroups of $G$ before getting to the incidence system.
To be more precise, let $I$ be a finite indexing set and $S = \{s_i \in G \mid i\in I\}$ such that $G = \langle S \rangle$ (that is, $S$ is a generating set for $G$).
For any $J\subseteq I$, define subgroups $G^{(J)} := \langle S^J \rangle$, where $S^J := \{s_j \mid j\in J\}$ and $G^{(\emptyset)} = \{e\}$.
\begin{definition}
    We say the pair $(G,S)$ satisfies the \defn{intersection property} if for any $J,K\subseteq I$ we have that     
    \begin{equation} \tag{IP}
    \label{eq:IP}
       G^{(J)}\cap G^{(K)} = G^{(J\cap K)}. 
    \end{equation}
    If the set $S$ of generators is clear from context, we will sometimes simply say that $G$ satisfies the intersection property. 
\end{definition}

Checking Equation~\eqref{eq:IP} for all the possible subsets $J,K \subset I$ is a hard task.
The following proposition shows how to do this more efficiently, provided we already know that the groups $G^{(I \setminus \{i\})}$ satisfies~\eqref{eq:IP}.

\begin{prop}\label{prop:IP_on_the_maximals}
    Consider a group $G$ with a finite generating set $S$ with index set $I$.
    Suppose that $(G^{(I\setminus\{i\})},S^{I\setminus\{i\}})$ satisfies~\eqref{eq:IP} for all $i\in\{1,\ldots,n\}$.
    Then $(G,S)$ satisfies~\eqref{eq:IP} if and only if $G^{(I\setminus\{i\})}\cap G^{(I\setminus\{j\})}= G^{(I\setminus\{i,j\})}$ for all $1 \leq  i, j \leq n.$
\end{prop}

\begin{proof}
    The forward implication is trivial, as the definition of intersection property covers the desired intersection.
    
    We thus only need to prove the backward implication. As all $(G^{(I\setminus\{i\})},S^{I\setminus\{i\}})$ satisfy the intersection property and $G^{(I\setminus\{i\})}\cap G^{(I\setminus\{j\})}= G^{(I\setminus\{i,j\})}$ for all $1 \leq  i, j \leq n$, most intersection are already covered. We nonetheless need to verify that $G^{(I\setminus\{i\})} \cap G^{(J)} = G^{(I\setminus\{i\}\cap J)}$, for any $J\subseteq I$. This is equivalent to proving that for any $J\subseteq I$, $G^{(I\setminus\{i\})} \cap G^{(I\setminus J)} = G^{(I\setminus [J\cup \{i\}])}$ and we can
    assume that $|J|>1$.
    Pick $j\in J$ with $j \neq i$. Then, notice that $G^{(I\setminus J)}\leq G^{(I\setminus\{j\})}$. Hence, $$G^{(I\setminus\{i\})} \cap G^{(I\setminus J)} = G^{(I\setminus\{i\})} \cap G^{(I\setminus\{j\})} \cap G^{(I\setminus J)} = G^{(I\setminus\{i,j\})} \cap G^{(I\setminus J)}.$$
    Since $(G^{(I\setminus\{j\})},S^{I\setminus\{j\}})$ satisfies the intersection property, we can conclude that $G^{(I\setminus\{i,j\})} \cap G^{(I\setminus J)} = G^{(I\setminus [J\cup \{i,j\}])}= G^{(I\setminus [J\cup \{i\}])}$.
\end{proof}

In some situations, the intersection property can be preserved under quotients.

\begin{lemma}\label{lem:general_IP_quotient_subs}
    Let $G$ be a group and let $A,B\leq G$ and $N\trianglelefteq G$.
    If $A\cap N = \{e\}=B\cap N$, then $ (A\cap B)/N\iso A\cap B \iso A/N \cap B/N$.
\end{lemma}
\begin{proof}
    The hypothesis guarantees that the quotient map $\pi \colon G \to G/N$ is an isomorphism when restricted to $A,B$ or $A \cap B$. We thus get that $A \cap B \iso \pi (A \cap B) = (A \cap B)/N$ but also that $\pi (A \cap B) = \pi(A) \cap \pi(B) = A/N \cap B/N$.
\end{proof}

\begin{prop}\label{prop:IP_for_quotients}
    Let $G$ be a group with a finite generating set $S = \{s_i \mid i\in I\}$ and $N\trianglelefteq G$ be a normal subgroup of $G$. 
    \begin{enumerate}
        \item \label{prop:IP:q1} Suppose $N\cap G^{(J)}= \{e\} =N\cap G^{(K)}$ for some $J,K\subset I$.
        If $G^{(J)}\cap G^{(K)} = G^{(J\cap K)}$, then $G^{(J)}/N\cap G^{(K)}/N \iso G^{(J\cap K)}/N$.
        \item \label{prop:IP:q2} Moreover, for any $J \subseteq I$, if $(G^{(J)},S^{J})$ satisfies~\eqref{eq:IP} and $G^{(J\setminus\{i\})}\cap N = \{e\}$ for every $i \in J$, then $(G^{(J)}/N,S^{J}/N)$ satisfies~\eqref{eq:IP}, where $S^{J}/N := \{s_iN\mid i\in J\}$.%\travis{I generalized this slightly for the statement that follows.}
    \end{enumerate}
\end{prop}

\begin{proof}
    Part (a) is a direct consequence of Lemma~\ref{lem:general_IP_quotient_subs}, with $A = G^{(J)}$ and $B = G^{(K)}$.
    
    Now we consider (b).
    Notice that, since $G^{J\setminus\{i\}} \cap N = \{e\}$ for every $i\in J$, we have $G^{(J\setminus\{i\})}/N\iso G^{(J\setminus\{i\})}$ for every $i\in J$.
    Additionally, since $(G^{(J)},S^{J})$ satisfies the intersection property~\eqref{eq:IP}, we see that the pairs $(G^{(J\setminus\{i\})},S^{J\setminus\{i\}})$ and $(G^{(J\setminus\{i\})}/N,S^{J\setminus\{i\}}/N)$ also satisfies~\eqref{eq:IP}. 
    Hence, by Proposition~\ref{prop:IP_on_the_maximals}, we need only to prove the intersection property in $G^{(J)}/N$ of the form $G^{(J\setminus\{i\})}/N \cap G^{(J \setminus \{j\})}/N$ for $i \neq j \in J$.
    That follows directly from~(a).
\end{proof}

We can reduce checking the intersection property~\eqref{eq:IP} to groups that cannot be written as a direct product.

\begin{coro}
    \label{coro:IP_products}
    Consider pairs $(G_1, S_1)$ and $(G_2, S_2)$.
    Both $(G_1, S_1)$ and $(G_2, S_2)$ satisfy the intersection property if and only if $(G_1 \times G_2, S_1 \sqcup S_2)$ satisfies the intersection property.
\end{coro}

\begin{proof}
    Let $I_1$ and $I_2$ be the index sets for $S_1$ and $S_2$, respectively.
    Let $G_{12} = G_1 \times G_2$.
    The forward direction is easy to see as $G_{12}^{(J)} \cap G_{12}^{(K)} = \{e\}$ whenever $J \subseteq I_1$ and $K \subseteq I_2$.
    The reverse direction follows from Proposition~\ref{prop:IP_for_quotients}(\ref{prop:IP:q2}) by taking $J = I_1$ (resp.\ $J = I_2$) and $N = G_2$ (resp.\ $N = G_1$).
\end{proof}

Now let us construct the specific coset incidence system associated to the pair $(G, S)$.
Let $\Sigma(G, S)$ denote the coset incidence system $\Sigma(G, (G_i)_{i \in I})$ with maximal parabolic subgroups
\[
G_i := \langle s_j \mid j\in I\setminus\{i\}\rangle = \langle S^{I\setminus \{i\}}\rangle = G^{(I\setminus\{i\})}.
\]
Note that in general $G^{(J)} = \langle S^J\rangle$ is distinct from $G^J := \langle G^j\mid j\in J\rangle$.
However, if $(G,S)$ satisfies the intersection property, both definitions coincide with $G^{(J)} = G^J$.
Indeed, by~\eqref{eq:IP}, it must then be that $G^i = \langle s_i\rangle$ and, for all $J\subseteq I$, that $G^J = \langle G^i\mid i \in J\rangle = \langle S^J\rangle = G^{(J)}$.

\subsection{Buekenhout Diagram}
Francis Buekenhout introduced in~\cite{buek} a diagram associated to incidence geometries. His insight was that, in many cases, the rank two residues of a geometry encode a lot of the global structure of the geometry. We will define this diagram for flag-transitive geometries, as this ensures that all residues of type $\{i,j\}$, for $i,j\in I$, are isomorphic to each other. 

Let $\Sigma$ be rank two geometry over $\{i,j\}$. The \defn{gonality} of $\Sigma$ is defined to be half of the girth of its incidence graph.
The \defn{$i$-diameter} of $\Sigma$ is the maximal distance in the incidence graph from $x$ to $y$, where $x$ is a element of type $i$ and $y$ is any element of $\Sigma$.
%The \defn{$j$-diameter} is defined in the same way, asking that $x$ is of type $j$ instead.
%\travis{I commented out the $j$-diameter as it wasn't saying anything new.}
\begin{definition}\label{def:BuekenhoutDiag}
   A \defn{Buekenhout diagram} $D$ is a complete graph whose vertex set is $I$ and
   \begin{itemize}
       \item each vertex $i \in I$ has a positive integer label $o_i$, called the \defn{i-order};
       \item each oriented edge $\{i,j\}$ has label a triple $(d_i,g_{i,j},d_j)$ of positive integers.
   \end{itemize}

   A incidence geometry $\Sigma$ is said to \defn{have diagram} $D$ if
\begin{itemize}
    \item for any flag $F$ of type $I\setminus\{i\}$, we have $|\Sigma_F| -1 = o_i$;
    \item all residues $\Sigma_F$ of type $\{i,j\}$ satisfy:
    \begin{enumerate}
        \item the gonality of $\Sigma_F$ is $g_{ij}$,
        \item the $i$-diameter (resp.~$j$-diameter) of $\Sigma_F$ is $d_i$ (resp.~$d_j)$.
    \end{enumerate} 
    
\end{itemize}  
\end{definition}

Equivalently, each rank 2 residue of type $\{i,j\}$ of $\Sigma$ must have a Buekenhout diagram
\[
\xymatrix@-1pc{*{\circ}\ar@{-}[rrr]^{d_{i},g_{ij},d_{j}}_(0.01){o_i}_(0.99){o_j} &&& *{\circ} }.
\]

Whenever $g_{ij}=d_{i}=d_{j}=n$, we say that this rank $2$ residue is a \defn{generalised $n$-gon} and, in this case, we omit $d_{i}$ and $d_{j}$. Moreover, if $g_{ij}=d_{i}=d_{j}=2$, we will call these residues a \defn{generalised digon}. The usual convention is to not draw at all the edges corresponding to generalized digons.

If $\Sigma$ is a coset geometry, we can determine the parameters $d_{i}$, $g_{ij}$, $d_{j}$ and $o_i$ using properties of the groups.
Let $f\in \mathbb{N}_0$ and $G^{i,j}_f$ denote product of $f$ subgroups, alternating between $G^i$ and $G^j$, more precisely $G^{i,j}_f := G^iG^jG^i\ldots$, with $f$ factors. For example, we have that $G^{i,j}_3 = G^iG^jG^i$ and $G^{j,i}_2=G^jG^i$. By definition, we set $G^{i,j}_0=G^{j,i}_0=B$.
The following results are gathered from~\cite{Pasini1994}.

\begin{prop}[{\cite[pp.~307]{Pasini1994}}] \label{prop:Pasini}
    Let $\Sigma = \Sigma(G,(G_0,G_1))$ be a rank $2$ firm and residually connected coset geometry.
    \begin{enumerate}
        \item The geometry $\Sigma$ is a generalized digon if and only if $G = G_0G_1=G_1G_0$.
        \item The point-diameter $d_{0}$ is equal to $f$ if and only if $G^{0,1}_f=G$ and $G^{0,1}_{f-1}\neq G$. Similarly, the line-diameter $d_1$ is equal to $f$ if and only if $G^{1,0}_f=G$ and $G^{1,0}_{f-1}\neq G$.
        \item The gonality $g = g_{01}$ is the minimal positive integer $f$ such that $G^{0,1}_f\cap G^{1,0}_f \neq G^{0,1}_{f-1}\cup G^{1,0}_{f-1}$.
        \item The $i$-orders are given by $o_i=[G_i:B]-1$, for $i\in \{0,1\}$.
    \end{enumerate}
\end{prop}

With the results above, we can describe all rank $2$ residues of any firm and residually connected coset geometry of rank $n\geq 2$ and describe its Buekenhout diagram.

\subsection{Shephard groups}\label{sec:back:subsec:ShephardGroup}

A simple graph $\Gamma$ is a set of vertices $V(\Gamma)$ and edges
\[
E(\Gamma) \subseteq \{\{x,y\} \mid x, y \in V(\Gamma) \text{ such that } x \neq y\};
\]
in particular, there is at most one edge between every pair of vertices and no loops.
Let $\Gamma$ be a simple graph such that every vertex is labeled by $p_i\in \mathbb{N}_{\geq 2}\cup \{\infty\}$ and every edge is labeled by $r_{i,j}\in \mathbb{N}_{\geq 3}\cup \{\infty\}$.
Moreover, whenever $r_{i,j}$ is odd, we require that $p_i=p_j$.
We call such a graph $\Gamma$ a \defn{Shephard presentation graph}.

For a presentation graph $\Gamma$, the corresponding \defn{Shephard group} $\Sh_{\Gamma}$ is the group given by the presentation
\[
\Sh_\Gamma = 
\left\langle
\begin{matrix}
    \\
 a_i, i\in V(\Gamma)   \\
    \\
\end{matrix}
\right\rvert
\left. \begin{matrix}
a_i^{p_i}= 1\textnormal{ if }p_i\neq \infty\\
\Delta_{a_i,a_j}^{r_{i,j}}=\Delta_{a_j,a_i}^{r_{i,j}}\textnormal{ for }\{i,j\}\in E(\Gamma)\textnormal{ and }r_{i,j}\neq \infty\\
\Delta_{a_i,a_j}^2=\Delta_{a_j,a_i}^2\textnormal{ if }\{i,j\}\notin E(\Gamma)
\end{matrix}
\right\rangle.
\]

\begin{remark}\label{remark:2labels}
   Each $\{i,j\} \notin E(\Gamma)$ corresponds to the relation $\Delta_{a_{i,j}}^2=\Delta_{a_{j,i}}^2$.
   Therefore, we could alternatively define $\Gamma$ as a complete graph with labels in $\mathbb{N}_{\geq 2} \sqcup \{\infty\}$, with the label $2$ corresponding to commuting generators.
   Hence, we will occasionally identify non-edges of $\Gamma$ with ``edges labeled $2$'' in $\Gamma$. 
\end{remark} 

Note that $\Sh_{\Gamma}$ has one generator for each vertex of $\Gamma$ and all relations involve at most one pair of generators. 
A Shephard group $\Sh_\Gamma$ is said to be of \defn{rank} $n$ if $\abs{V(\Gamma)} = n$.
We say $\Sh_{\Gamma}$ is \defn{irreducible} if $\Gamma$ is connected and otherwise it is \defn{reducible}.
If $\Gamma$ has connected components $\Gamma_1, \dotsc, \Gamma_N$, then clearly
$
\Sh_{\Gamma} \iso \prod_{i=1}^N \Sh_{\Gamma_i}.
$
We will denote rank $2$ Shephard groups by $\Sh_\Gamma = \Sh(p_1,r,p_2)$, where the presentation graph $\Gamma$ is 
\[
\xymatrix@-1pc{*{\bullet}\ar@{-}[rr]^{r}^(0.01){p_1}^(0.99){p_2} && *{\bullet}}
%\qquad (r > 2),
%\hspace{60pt}
%\xymatrix@-1pc{*{\bullet}\ar@{}[rr]^(0.01){p_1}^(0.99){p_2} && *{\bullet}}
%\qquad (r = 2).
\]
where by Remark~\ref{remark:2labels}, when $r = 2$ this could also be draw as two disjoint vertices.
Note that $\Sh(p_1,r,p_2)$ is irreducible if and only if $r > 2$.
Furthermore, clearly $\Sh(p_1,r,p_2) \iso \Sh(p_2,r,p_1)$, and so without loss of generality, we could assume $p_1 \leq p_2$.

Shephard groups generalize both Coxeter groups and Artin--Tits groups.
Indeed, if $p_i=2$ for all $i\in V(\Gamma)$, then $\Sh_{\Gamma}$ is a Coxeter group, while if all $p_i=\infty$, then $\Sh_\Gamma$ is an Artin--Tits group.
Furthermore, every Shephard group is a quotient of the corresponding Artin--Tits group.
Finite Shephard groups are given by the classification of complex reflection groups~\cite{ShephardTodd1954}, which were then studied by Coxeter in~\cite{Coxeter1975-vz} as the groups acting on finite complex polytopes.

Let us now explain why we imposed the condition that $p_i = p_j$ whenever $r_{ij}$ is odd. The following appears to be well-known to experts, being mentioned in~\cite[Section 8.9]{Coxeter1975-vz}.

\begin{lemma}
    \label{lem:trivial_group}
    If $r$ is odd and $p_1 \neq p_2$, let $g = \gcd(p_1, p_2)$ then $\Sh(p_1,r,p_2) = \Sh(g,r,g)$.
    In particular, if $g = 1$, then $\Sh(p_1,r,p_2) = \{e\}$ is the trivial group.
\end{lemma}

\begin{proof}
    Without loss of generality assume $p_1 < p_2$, and we have
    \[
    \Delta_{a_1,a_2}^r = a_1^{p_1} \Delta_{a_1,a_2}^r = a_1^{p_1-1} \Delta_{a_1,a_2}^r a_2 = \cdots = \Delta_{a_1,a_2}^r a_2^{p_1}
    \qquad \Longrightarrow \qquad
    a_2^{p_1} = 1.
    \]
    Hence $a_2^g = 1$, and we can apply the same argument swapping to show that $a_1^g = 1$.
\end{proof}

The Shephard group $\Sh_{\Gamma}$ has a distinguished set of generators $S = \{a_i \mid i \in V(\Gamma)\}$, and so we denote the corresponding coset incidence system by $\Sigma(\Sh_{\Gamma}) := \Sigma(\Sh_{\Gamma}, S)$.

In concurrence with the definition of right-angled Coxeter/Artin--Tits groups, we make the following definition.

\begin{definition}\label{def:RAShepG}
    Let $\Gamma$ be a Shephard presentation graph.
    The Shephard group $\Sh_\Gamma$ is said to be \defn{right-angled} if the labels of all edges of $\Gamma$ is $\infty$.
\end{definition}

%We say $\Lambda$ is an \defn{induced subgraph of $\Gamma$} if $V(\Lambda)\subseteq V(\Gamma)$ and, for every pair of vertices $\{i,j\}\in V(\Lambda)$, if $\{i,j\}\in E(\Gamma)$, then $\{i,j\}\in E(\Lambda)$, with the exact same labeling.
If $\Lambda$ is an induced subgraph of $\Gamma$, note that $\Lambda$ is also a Shephard presentation graph.
However, the subgroup of $\Sh_\Gamma$ generated by $\langle a_i\mid i\in \Lambda\rangle$ is not necessarily the same group as $\Sh_\Lambda$.
Hence, we will distinguish these two groups. Let $\Sh_\Gamma(\Lambda)$ be the subgroup of $\Sh_\Gamma$ generated by the vertices of $\Lambda$.
In other words, we have $\Sh_\Gamma(\Lambda):= \langle a_i \mid i\in \Lambda\rangle \subseteq \Sh_\Gamma$.

\begin{example}
Consider the group $\Sh(3,3,5) = 1$ (Lemma~\ref{lem:trivial_group}), and let $\Gamma$ be the corresponding graph.
Necessarily all $\Sh_{\Gamma}(\Lambda) = 1$ for any induced subgraph $\Lambda$, but when $V(\Lambda) = \{1\}$ (resp.\ $V(\Lambda) = \{2\}$), then we have $\Sh_{\Lambda} = \ZZ/3\ZZ$ (resp.\ $\Sh_{\Lambda} = \ZZ/5\ZZ$).
However, in this case, the graph $\Gamma$ is not a Shephard presentation graph.
\end{example}

From~\cite{Goldman2024,2ShepKatherine} and private communication with Katherine Goldman, we know that for some induced subdiagrams $\Lambda\subseteq \Gamma$ of particular Shephard groups $\Sh_\Gamma$, we have that $\Sh_\Gamma(\Lambda)=\Sh_\Lambda$.
However, it is not usually proven for every possible induced subdiagram.

\begin{definition}
    The Shephard group associated to a Shephard presentation graph $\Gamma$ is \defn{tame} if for any induced subgraph $\Lambda\subseteq \Gamma$, we have that $\Sh_{\Gamma}(\Lambda) = \Sh_{\Lambda}$.
    %\travis{Since you said this term was not defined previously.}
\end{definition}

To our understanding, it is conjectured that all Shephard groups are tame as no counterexample is known (when $\Gamma$ is a Shephard presentation graph).
Vacuously, any rank $1$ Shephard group is tame, where it is just a cyclic group and has no proper subgraphs.

\begin{prop}
    \label{prop:tame_product}
    Consider a Shephard presentation graph $\Gamma$ whose connected components are $\Gamma_1, \dotsc, \Gamma_N$.
    Then the Shephard group $\Sh_{\Gamma}$ is tame if and only if the Shephard groups $\Sh_{\Gamma_i}$, for all $i = 1, \dotsc, N$ are tame.
\end{prop}

\begin{proof}
    This follows immediately from the fact $\Sh_{\Gamma} \iso \Sh_{\Gamma_1} \times \cdots \times \Sh_{\Gamma_N}$ and basic properties of direct products of groups.
\end{proof}

Thus, it is sufficient to show that all irreducible Shephard groups are tame.

\section{Incidence Geometries for Rank Two Shephard Groups}

In this section, we investigate the natural rank two coset geometries associated to a Shephard group with two generators.
Hence, we consider the Shephard groups $\Sh_\Gamma = \Sh(p_1, r, p_2)$.
%\langle a_1,a_2\mid a_1^{p_1}=a_2^{p_2}=1, \Delta_{a_1,a_2}^{r}=\Delta_{a_2,a_1}^{r}\rangle$, with $p_1,p_2,r\in \mathbb{N}_{\geq 2}\cup \{\infty\}$, whose %Shephard presentation graph $\Gamma$ is 
%$$\xymatrix@-1pc{*{\bullet}\ar@{-}[rr]^{r}^(0.01){p_1}^(0.99){p_2} && *{\bullet}}.$$
Recall that, if $r$ is odd, we must have that $p_1 = p_2$.

Each rank $2$ Shephard group $\Sh(p_1,r,p_2)$ is a quotient of the Shephard group $\Sh(p_1,\infty,p_2)$. Note that $\Sh(p_1,\infty,p_2)$ is the free product of two cyclic groups $\Cyc_{p_1}*\Cyc_{p_2}$.
Every element $w \in \Sh(p_1,\infty,p_2)=\langle a_1,a_2\rangle$ has a unique normal form $w = a_{k_1}^{l_1}a_{k_2}^{l_2}\ldots a_{k_m}^{l_m}$ with $k_{i}\in\{1,2\}$, $k_{i}\neq k_{i+1}$ for $1\leq i\leq m-1$, and $1\leq l_i< p_{k_i}$ (for more details, see~\cite{Lyndon2001}).
The syllable length of the element $w$ above, denoted $\abs{w}$, is equal to $m$.%\travis{This is a highly nonstandard definition of length. I would call it syllable length.}
By the definition of free product of groups, it can be easily proved that $\Sh(p_1,\infty,p_2)$ is tame and satisfies the intersection property. We will now prove for the remaining cases.

\begin{lemma}\label{lem:LengthReducedWordsRank2}
    Let $\Sh(p_1, r, p_2)$ be a rank $2$ Shephard group with $r\neq \infty$.
    Let $w$ be a non-trivial element of $\Sh(p_1,\infty,p_2)$ and consider the natural projection $\pi\colon \Sh(p_1,\infty,p_2) \to \Sh(p_1,r,p_2)$.
    If $\pi(w) = e$ in $\Sh(p_1,r,p_2)$, then $\abs{w} \geq 2r$.
\end{lemma}

\begin{proof}
The case of $p_1,p_2 \neq \infty$ is covered by~\cite[Proposition~4.3]{2ShepKatherine}. If both $p_1,p_2=\infty$, then $\Sh_\Gamma$ is an Artin--Tits group and the result is given in~\cite[Lemma~6]{Appel1983}.
Suppose without loss of generality that $p_1=\infty$ and $p_2\neq \infty$.
Suppose $w$ is a non-trivial element of $\Sh(\infty,\infty,p_2)$ such that $\abs{w} < 2r$ and $\pi(w)=e$. 
Note that, as $p_1 \neq p_2$, we must have that $r$ is even (see Lemma~\ref{lem:trivial_group}).
Then, there is also a surjective homomorphism $\lambda \colon \Sh_\Gamma = \Sh(\infty,r,p_2) \to \Sh(q,r,p_2)$ for any $q \in \mathbb{N}_{\geq 2}$.
The image of $w$ under $\lambda\circ \pi$ is still representing the identity in $\Sh(q,r,p_2)$. Consider $q$ to be greater than the biggest exponent of $a_1$ appearing in $w$.
Then there are two surjective homomorphisms $\gamma \colon \Sh(\infty,\infty,p_2)\to \Sh(q,\infty,p_2)$ and $\eta \colon \Sh(q,\infty,p_2) \to \Sh(q,r,p_2)$.
Since $q$ is greater than any exponent of $a_1$ in $w$, we have that $\gamma(w) = w$.
Hence, from the case above of $p_1,p_2\neq \infty$, it must be that $\eta\circ\gamma(w)\neq e$, a contradiction with the fact that $\lambda\circ \pi(w)=e$. 
\end{proof}

\begin{lemma}\label{lem:Rank2_properly_defined}
    Let $G= \Sh(p_1,\infty,p_2) = \langle a_1,a_2\mid a_1^{p_1} = a_2^{p_2}=1\rangle$.
    We have $N \cap \langle f_1\rangle = \{e\}  = N \cap \langle f_2\rangle$, where $N\trianglelefteq G$ is the normal closure of the subgroup $\langle \Delta_{f_1,f_2}^r(\Delta_{f_2,f_1}^r)^{-1}\rangle$ for $r\neq \infty$.
\end{lemma}

\begin{proof}
    By contradiction, suppose that there exists an $i\in\{1,2\}$ such that $N\cap \langle a_i\rangle\neq \langle 1\rangle$. 
    Then, we have $a_i^{k}\in N$, for $0<k<p_i$. This would mean that $\pi(a_i^k)=e$, where $\pi$ is defined as in Lemma~\ref{lem:LengthReducedWordsRank2}.
    However $w=a_i^k$ is a non-trivial element of $G$, with $k$ not a multiple of $p_i$, whose length is equal to one and, by the intersection with $N$, $\pi(w)=e$.
    This contradicts Lemma~\ref{lem:LengthReducedWordsRank2}, and thus concludes the proof.
\end{proof}

Note that the condition forcing $p_1$ to be equal to $p_2$ whenever $r$ is odd is essential for the previous results to hold. 
Recall that a rank two Shephard group $\Sh_\Gamma = \Sh(p_1,r,p_2)$ satisfies the intersection property if $\langle a_1 \rangle \cap \langle a_2 \rangle = \{e\}$.

\begin{thm}\label{thm:Shephard_IP_rank2}
    Let $\Sh_\Gamma= \Sh(p_1,r,p_2)$ be a rank two Shephard group. Then $\Sh_\Gamma$ is tame and satisfies the intersection property.
\end{thm}

\begin{proof}
If $r=\infty$, $\Sh_\Gamma$ is a free product and the results is well known in this case.

Suppose now that $2\leq r<\infty$.
In this case, $\Sh_\Gamma=\Sh(p_1,r,p_2)$ is the quotient of $\Sh(p_1,\infty,p_2)$ by the normal closure $N=\langle\langle \Delta_{a_1,a_2}^r(\Delta_{a_2,a_1}^r)^{-1}\rangle\rangle$. 
By Lemma~\ref{lem:Rank2_properly_defined}, we have that $\langle a_1\rangle \cap N = \{ e \} = \langle a_2 \rangle \cap N$.
Hence, the order of $a_i$ after the quotient of $\Sh(p_1,\infty,p_2)$ by $N$ is still $p_i$. Therefore, $\Sh_\Gamma(\{i\}) = \langle a_i\rangle = \Sh_{\{i\}}$, so that $\Sh_\Gamma$ is tame for any rank $2$ $\Gamma$.
Following Proposition~\ref{prop:IP_for_quotients}, we  can conclude that $\Sh(p_1,r,p_2)$ satisfies the intersection property as well.
\end{proof}

\begin{comment}
We have thus proved that all rank $2$ Shephard groups satisfy the intersection property. Additionally, any rank two coset incidence system $(G,(G_i)_{i\in I})$ is flag-transitive (\cite[Lemma 1.8.6]{buekenhout2013diagram}). Therefore, by the intersection property, we have that it is also residual-connected.
Therefore we have the following Corollary of the above result.

\begin{coro}\label{coro:ShephardFTRC_CG_rank2}
    The coset incidence system of a rank $2$ Shephard group is a flag-transitive, residually connected firm coset geometry.
\end{coro}
\begin{proof}
    
    FIRMNESS COMES FROM THE FACT THAT $B=\langle 1\rangle$ AND $p,q\geq 2$, HENCE $|G^i:B|=|G^i| = |p|\ (or\ |q|)$!
\end{proof}

An immediate Corollary of Theorem~\ref{thm:Shephard_IP_rank2} is the following.

\begin{coro}\label{coro:InducedSubRank2}
    Let $\Lambda$ be an induced subdiagram of a rank 2 diagram $\Gamma$. Then $\Sh_\Gamma(\Lambda) \iso \Sh_\Lambda$.
\end{coro}
\end{comment}

Consider now the coset geometry $\Sigma(p_1,r,p_2) = \Sigma(\Sh(p_1,r,p_2))$. %\travis{Do we actually know this is a geometry? I think it follows from the intersection property, but that is never explicitly stated.}\claudio{In a way, it could follow from the intersection property. Nonetheless, any coset incidence system of rank $2$ is trivially a geometry.}
We will describe the rank 2 geometries $\GG$ obtained in terms of point and line diameters and gonality.

\begin{prop}\label{prop:gon}
    The gonality of $\Sigma(p_1,r,p_2)$ is $r$.
\end{prop}

\begin{proof}Let $G=\Sh(p_1,r,p_2)$.
    To prove that the gonality of $\Sigma(p_1,r,p_2)$ is $r$, we will use Proposition~\ref{prop:Pasini}. 
    First, let us consider the case of $f < r$. We wish to show that 

    \begin{equation}\label{eq:gonality}
        G^{1,2}_f\cap G^{2,1}_f = G^{1,2}_{f-1}\cup G^{2,1}_{f-1}.
    \end{equation}
    Notice that, given the tameness of the rank $2$ Shephard groups, the image by $\pi$ of an non-trivial element $w = a_{k_1}^{l_{k_1}}a_{k_2}^{l_{k_2}}\ldots a_{k_f}^{l_{k_f}}$ of $\Sh(p_1,\infty,p_2)$ must belong to $G_f^{k_1,k_2}$.
    
    By contradiction, assume that equality  of equation~\ref{eq:gonality} does not hold. Note that $G^{0,1}_f\cap G^{1,0}_f$ always contains $G^{0,1}_{f-1}\cup G^{1,0}_{f-1}$. Thus, we can find an element $g \in G^{0,1}_f\cap G^{1,0}_f $ which is not in $G^{0,1}_{f-1}\cup G^{1,0}_{f-1}$. This means that there must exist two non-trivial elements, $w,w'$ of $\Sh(p_1,\infty,p_2)$ of length $f$ such that $\pi(w)=\pi(w')=g$, one starting with $a_1$ and the other starting with $a_2$. Therefore, we have
    
    $$ g = \pi(a_1^{k_1}a_2^{k_2}a_1^{k_3} \cdots a_{i_f}^{k_f}) = \pi(a_2^{m_1}a_1^{m_2}a_2^{m_3} \cdots a_{j_f}^{m_f})$$
    where $i_f = 1, j_f = 2$ if $f$ is odd and $i_f =2, j_f =1$ if $f$ is even. In any case, the element $ww'^{-1}=(a_1^{k_1}a_2^{k_2}a_1^{k_3} \cdots a_{i_f}^{k_f})(a_2^{m_1}a_1^{m_2}a_2^{m_3} \cdots a_{j_f}^{m_f})^{-1}$ is a non-trivial element of length $2f$ of $\Sh(p_1,\infty,p_2)$ for which $\pi(ww'^{-1})$ is the identity. Since $f < r$, this contradicts Lemma \ref{lem:LengthReducedWordsRank2}. Therefore, it must be that $G^{0,1}_f\cap G^{1,0}_f = G^{0,1}_{f-1}\cup G^{1,0}_{f-1}$ for any $f < r$.

    Now let $f = r$. We need to show the equality in Equation~\eqref{eq:gonality} does not hold for $f = r$. We consider the words $\Delta^r_{a_1,a_2}$ and $\Delta^r_{a_2,a_1}$ used as the defining relation of the Shephard group. These words have length $r$ and thus the element $g = \Delta^r_{a_1,a_2} = \Delta^r_{a_2,a_1}$ of $\Sh(p_1,r,p_2)$ belongs to $G^{1,2}_f\cap G^{2,1}_f$. Suppose that $g=\Delta^r_{a_1,a_2} = \Delta^r_{a_2,a_1}$ also belongs to $G^{1,2}_{f-1}\cup G^{2,1}_{f-1}$. It means that there exists $w\in \Sh(p_1,\infty,p_2)$ of length $f-1$ such that $\pi(w)=g$. In $\Sh(p_1,\infty,p_2)$, one of the words $\Delta^r_{a_1,a_2} w^{-1}$ or $\Delta^r_{a_2,a_1} w^{-1}$ must be in the normal form and has length strictly smaller than $2r$. Since the projection of this word into $\Sh(p_1,r,p_2)$ gives the identity, this again contradicts Lemma \ref{lem:LengthReducedWordsRank2}, and concludes the proof.    
\end{proof}

\begin{prop}\label{prop:r=2digon}
    If $r=2$, then $\Sigma(p_1,r,p_2)$ is a generalised digon.
\end{prop}

\begin{proof}
    Let $\Sh(p_1,2,p_2) = \langle a_1^{p_1} = a_2^{p_2}= 1, \Delta_{a_1,a_2}^2=\Delta_{a_2,a_1}^2\rangle$.
    Notice that the Artin--Tits relation implies that $a_1$ and $a_2$ commute.
    Hence, $\Sh(p_1,2,p_2) \iso \langle a_1\rangle\times \langle a_2 \rangle$.
    We can then immediately conclude using the first item of Proposition \ref{prop:Pasini}.
\end{proof}

\begin{prop}
     Let $\Sh(p_1,r,p_2)$ be a finite Shephard group of rank two.
     Then $\Sigma(p_1,r,p_2)$ is one of the geometries listed in Table~\ref{table:finiteR2}.
\end{prop}

\begin{table}[t]
\begin{tabular}{|c|c|cc}
\hline
  Shephard group; ST &
  Buekenhout diagram &
  \multicolumn{1}{c|}{Shephard group; ST} &
  \multicolumn{1}{c|}{\text{Buekenhout diagram}} \\ \hline
$\Sh(2,r,2)$; $G(r,r,2)$  & $\xymatrix@-1pc{*{\circ}\ar@{-}[rr]^{r,r,r}_(0.01){1}_(0.99){1} && *{\circ}}$   & \multicolumn{1}{c|}{$\Sh(3,3,3)$; $G_4$} & \multicolumn{1}{c|}{$\xymatrix@-1pc{*{\circ}\ar@{-}[rr]^{4,3,4}_(0.01){2}_(0.99){2} && *{\circ}}$}   \\
\hline
$\Sh(2,4,p)$; $G(p,1,2)$  & $\xymatrix@-1pc{*{\circ}\ar@{-}[rr]^{4,4,4}_(0.01){1}_(0.99){p-1} && *{\circ}}$    & \multicolumn{1}{c|}{$\Sh(3,4,3)$; $G_5$} & \multicolumn{1}{c|}{$\xymatrix@-1pc{*{\circ}\ar@{-}[rr]^{6,4,6}_(0.01){2}_(0.99){2} && *{\circ}}$}   \\ \hline
$\Sh(2,6,3)$; $G_6$  & $\xymatrix@-1pc{*{\circ}\ar@{-}[rr]^{8,6,8}_(0.01){1}_(0.99){2} && *{\circ}}$    & \multicolumn{1}{c|}{$\Sh(3,4,4)$; $G_{10}$} & \multicolumn{1}{c|}{$\xymatrix@-1pc{*{\circ}\ar@{-}[rr]^{8,4,8}_(0.01){2}_(0.99){3} && *{\circ}}$}   \\ \hline
$\Sh(2,6,4)$; $G_9$  & $\xymatrix@-1pc{*{\circ}\ar@{-}[rr]^{10,6,10}_(0.01){1}_(0.99){3} && *{\circ}}$  & \multicolumn{1}{c|}{$\Sh(3,5,3)$; $G_{20}$} & \multicolumn{1}{c|}{$\xymatrix@-1pc{*{\circ}\ar@{-}[rr]^{11,5,11}_(0.01){2}_(0.99){2} && *{\circ}}$} \\ \hline
$\Sh(2,6,5)$; $G_{17}$  & $\xymatrix@-1pc{*{\circ}\ar@{-}[rr]^{14,6,14}_(0.01){1}_(0.99){4} && *{\circ}}$  & \multicolumn{1}{c|}{$\Sh(3,4,5)$; $G_{18}$} & \multicolumn{1}{c|}{$\xymatrix@-1pc{*{\circ}\ar@{-}[rr]^{12,4,12}_(0.01){2}_(0.99){4} && *{\circ}}$} \\ \hline
$\Sh(2,8,3)$; $G_{14}$ & $\xymatrix@-1pc{*{\circ}\ar@{-}[rr]^{12,8,12}_(0.01){1}_(0.99){2} && *{\circ}}$ & \multicolumn{1}{c|}{$\Sh(4,3,4)$; $G_8$}             & \multicolumn{1}{c|}{$\xymatrix@-1pc{*{\circ}\ar@{-}[rr]^{5,3,5}_(0.01){4}_(0.99){4} && *{\circ}}$}   \\ \hline
$\Sh(2,10,3)$; $G_{21}$ & $\xymatrix@-1pc{*{\circ}\ar@{-}[rr]^{22,10,22}_(0.01){1}_(0.99){2} && *{\circ}}$ &  \multicolumn{1}{c|}{$\Sh(5,3,5)$; $G_{16}$} & \multicolumn{1}{c|}{$\xymatrix@-1pc{*{\circ}\ar@{-}[rr]^{7,3,7}_(0.01){4}_(0.99){4} && *{\circ}}$}        \\ \hline
$\Sh(p_1,2,p_2)$  & $\xymatrix@-1pc{*{\circ}\ar@{-}[rr]^{2,2,2}_(0.01){p_1-1}_(0.99){p_2-1} && *{\circ}}$    &  &        \\ \cline{1-2}
\end{tabular}
\smallskip
\caption{Buekenhout diagram of finite Shephard groups of rank two, where we have also included the Shephard--Todd (ST) classification of the groups, for any $2 \leq p, p_1, p_2, r < \infty$.}
\label{table:finiteR2}
\end{table}

\begin{proof}
Recall that the list of finite Shephard groups is well known. We therefore only need to find the parameters of the rank two coset geometries associated to each of them.
Proposition~\ref{prop:gon} gives us the gonality of all these geometries. We now investigate the points and line diameters.

 First, note that $\Sh(2,r,2)$ is a Coxeter group, and thus the associated coset geometry is an $r$-gon, and thus the point and line diameter are also equal to $r$. Additionally, by Proposition~\ref{prop:r=2digon}, the associated geometry of $\Sh(p_1,2,p_2)$ is a generalized digon.

 Consider now the group $\Sh(2,4,p):=\langle a_1,a_2\mid a_1^2=a_2^p=1, \Delta_{a_1,a_2}^4=\Delta_{a_2,a_1}^4\rangle$.
 As the gonality is $4$, we have that the diameters must be at least $4$.
 From the Artin--Tits relation, we have that $a_1a_2a_1a_2 = a_2a_1a_2a_1$.
 Note that this implies that $a_2$ commutes with $a_1a_2a_1$.
 As $(a_1a_2a_1)^k = a_1a_2^ka_1$, we have also that $a_2^la_1a_2^ka_1 = a_1a_2^ka_1a_2^l$, for $1 \leq k,l < p$.
 
Hence, any reduced word of $\Sh(2,4,p)$ has at most length $4$.
Thus, by Proposition~\ref{prop:Pasini}(b), we conclude that both the point and line diameters are at most $4$.
They must then be equal to $4$ as they must be at least as big as the gonality.
Note that this argument holds also when $p = \infty$ so that $\Sigma(\infty,4,2)$ is also a generalized $4$-gon, but is infinite.

The diameters for the remaining cases have been computationally determined by using either Magma or GAP.
\end{proof}

\section{Higher Ranks}\label{sec:higher}

\subsection{Tame Shephard Groups}

In this section, we investigate whether Shephard groups satisfy the intersection property. We will show that they do if we assume that their parabolic subgroups are exactly what one would expect.

\begin{thm}\label{thm:Tame_imply_IP}

If $\Gamma$ is tame, then $\Sh_\Gamma$ satisfies the intersection property.
\end{thm}
\begin{proof}
    Let $r$ be the rank of $\Gamma$. We will proceed by induction on the rank. There is nothing to prove for rank 1, and the case rank 2 is handled by Theorem~\ref{thm:Shephard_IP_rank2}.

    We now suppose that the theorem holds for any $\Gamma$ of rank strictly less than $r$. Let $\Lambda_1$ and $\Lambda_2$ be any two proper sub-diagrams of $\Lambda$ and let $\Gamma'$ be the sub-diagram induced by $\Lambda_1 \cup \Lambda_2$. If $\Gamma'$ is not equal to $\Gamma$, we are done by applying the induction hypothesis. We can thus suppose that $\Gamma = \Gamma'$. The groups $\Sh_\Gamma(\Lambda_1)$ and $\Sh_\Gamma(\Lambda_2)$ share a common subgroup isomorphic to $\Sh_\Gamma(\Lambda_1 \cap \Lambda_2)$. We can therefore construct the free-amalgamated product $G = \Sh_\Gamma(\Lambda_1) *_{\Sh_\Gamma(\Lambda_1 \cap \Lambda_2)} \Sh_\Gamma(\Lambda_2)$. By the tameness of $\Gamma$, we have that $\Sh_\Gamma(\Lambda_i) \iso \Sh_{\Lambda_i}$ for $i = 1,2$ and that $\Sh_\Gamma(\Lambda_1 \cap \Lambda_2) \iso \Sh_{\Lambda_1 \cap \Lambda_2}$. Therefore, the group $G$ is also a Shephard group whose diagram is the join of $\Lambda_1$ and $\Lambda_2$, with infinity-labeled edges between generators of $\Lambda_1$ and $\Lambda_2$.a
    
    In $G$, we have that $G_1 \cap G_2 = G_{1,2} $ where $G_i$ is the canonical image of $\Sh_\Gamma(\Lambda_i)$ in $G$ and $G_{1,2}$ is the image of $\Sh_\Gamma(\Lambda_1 \cap \Lambda_2)$. By Lemma \ref{lem:general_IP_quotient_subs}, we can conclude that the same holds in $\Sh_\Gamma$, as $\Sh_\Gamma$ is the quotient of $G$ by the Shephard relations between generators in $\Lambda_1$ and $\Lambda_2$. Therefore, we have that $\Sh_\Gamma(\Lambda_1) \cap \Sh_\Gamma(\Lambda_2) = \Sh_\Gamma(\Lambda_1 \cap \Lambda_2)$, as desired.
\end{proof}

\begin{thm}\label{thm:Shephard_finite_tame}
    If $\Sh_\Gamma$ is a finite group, then $\Sh_\Gamma$ is tame.
\end{thm}
\begin{proof}
    If $\Sh_\Gamma$ is a finite Coxeter group, this is done. If not, this is proven by Coxeter in~\cite{Coxeter1975-vz}, and restated by Goldman in~\cite{Goldman2024}.
\end{proof}

\subsection{Rank 3}

Consider now a rank 3 Shephard group with the Shephard presentation graph $\Gamma$ given by
\[
\xymatrix@-0.5pc{*{\bullet}\ar@{-}[rr]^{r_{1,2}}^(0.01){p_1}^(0.99){p_2}&& *{\bullet}
\\
 & *{\bullet}\ar@{-}[ru]_{r_{2,3}}_(0.01){p_3}\ar@{-}[lu]^{r_{1,3}}}
\]
%\travis{Here's an issue where not allowing labels of $2$ in the graph is a problem. We probably should equate label 2 edges with non-edges.}\claudio{I will add in the description of the Shephard representation graph that label $2$ is sort of related with non-edges.}
and hence we consider the group $\Sh_\Gamma := \langle a_i , i\in\{1,2,3\}\mid a_i^{p_i} = 1, \Delta_{a_i,a_j}^{r_{i,j}}=\Delta_{a_j,a_i}^{r_{i,j}}\rangle$, where we allow some of the $r_{i,j}$ to be equal to $2$ (see Remark~\ref{remark:2labels}).
Note that, contrary to Coxeter and Artin--Tits groups, we still have yet to prove the tameness.
The following result gives conditions for which we can guarantee tameness and, by Theorem~\ref{thm:Tame_imply_IP}, the intersection property.

\begin{thm}
    Let $\Gamma$ be a rank $3$ diagram such that $\Sh_{\Gamma}$ is an irreducible Shephard group.
    Then $\Sh_\Gamma$ is tame if one of the following is true:
    \begin{enumerate}
        \item $\Sh_\Gamma$ is a finite group.
        \item For all induced subgraphs $\Lambda$ of $\Gamma$, with $\Lambda \neq \Gamma$, $\Sh_\Lambda$ is a finite group.
        \item $\frac{1}{{r_{1,2}}} + \frac{1}{{r_{2,3}}} + \frac{1}{{r_{1,3}}} \leq 1$, if $r_{i,j}\neq \infty$ for all $(i,j)\in\{(1,2),(2,3),(1,3)\}$.
        \item There is $r_{i,j}=\infty$, for $(i,j)\in\{(1,2),(2,3),(1,3)\}$.
        %\item $(r_{1,2},r_{2,3},r_{1,3}) \in \{(2,2,k), (2,k,2), (k,2,2)\}$, for $k\in \mathbb{N}_{\geq 2}\cup \{\infty\}$.
        %\claudio{should we change this and write instead "$\Sh_\Gamma$ is reducible"?}\travis{I think we should make it irreducible and have changed it accordingly. I added a note about the reducible case below.}\claudio{agreed with the change.}
    \end{enumerate}
\end{thm}

\begin{proof}
    We prove the result case-by-case.

    \smallskip \noindent \textbf{Case (a):}
    This is Theorem~\ref{thm:Shephard_finite_tame}.

    \smallskip \noindent \textbf{Case (b):}
    This is a consequence of~\cite[Theorem~A]{Goldman2024}.
    
    \smallskip \noindent \textbf{Case (c):}
    The case where all $p_i$ are finite is covered by~\cite[Theorem~B]{2ShepKatherine}.
    The proof presented here covers all cases, independently of the values of the $p_i$.
    Let $G = \Sh_\Gamma$ such that all $r_{i,j}\neq \infty$ for all $(i,j)\in\{(1,2),(2,3),(1,3)\}$. Consider $X$, $Y$ and $Z$ to be the following rank 2 Shephard groups:
    $X = \langle a_1,a_2\mid a_1^{p_1}=a_2^{p_2}=1, \Delta_{a_1,a_2}^{r_{1,2}}=\Delta_{a_2,a_1}^{r_{1,2}}\rangle$,
    $Y = \langle a_2,a_3\mid a_2^{p_2}=a_3^{p_3}=1, \Delta_{a_2,a_3}^{r_{2,3}}=\Delta_{a_3,a_2}^{r_{2,3}}\rangle$ and $Z = \langle a_1,a_3\mid a_1^{p_1}=a_3^{p_3}=1, \Delta_{a_1,a_3}^{r_{1,3}}=\Delta_{a_3,a_1}^{r_{1,3}}\rangle$.
    Note that we are not considering the subgroups of $G$ generated by the above generators, but actually taking distinct groups with the presentation above. Hence, we can define the following triangle of groups and injective homomorphisms:
    $$
    \xymatrix@-1pc{&& X && \\
    \\
    &A\ar@{->}[ruu]\ar@{->}[ldd]& *{\langle 1\rangle}
    \ar@{->}[r] 
    \ar@{->}[l] 
    \ar@{->}[dd]
    & B \ar@{->}[luu]\ar@{->}[rdd]&\\
    \\
    Z && C\ar@{->}[rr]\ar@{->}[ll] && Y
    }
    $$
    where $A=\langle a_1\mid a_1^{p_1}=1\rangle$, $B=\langle a_2\mid a_2^{p_2}=1\rangle$ and $C=\langle a_3\mid a_3^{p_3}=1\rangle$. As all rank 2 diagrams are tame, we get that $X$, $Y$ and $Z$ satisfy the intersection property and each subgroup $\langle a_i\mid a_i^{p_i}\rangle$ embeds naturally into each of these three groups.
    Hence, by the results given in rank 2, we have that $A$, $B$, $C$ inject naturally into $X$, $Y$ and $Z$ (considering the triangle of groups above).
    Let $G$ be the colimit of this triangle of groups, which by~\cite{stallings_1991}, is the Shephard group defined by $\Gamma$.
    
    By~\cite{stallings_1991}, we have that the natural maps of $X\to G$, $Y\to G$ and $Z\to G$ are injective if the triangle of groups if non-positively curved. Following~\cite{stallings_1991}, this property can be translated as 
    $\theta_{1,2}+\theta_{2,3}+\theta_{1,3}\leq \pi$, where $\theta_{i,j}$ is defined as the ``angle'' between the groups $\langle a_i\mid a_i^{p_i} = 1\rangle$ and $\langle a_j\mid a_j^{p_i} = 1\rangle$ in $\langle a_i, a_j\mid a_i^{p_i}=a_j^{p_j}=1, \Delta_{a_i,a_j}^{r_{i,j}}=\Delta_{a_j,a_i}^{r_{i,j}}\rangle$.
    To compute this angle, consider
    \[
    \phi_{i,j} \colon \langle a_i\mid a_i^{p_i} = 1\rangle*\langle a_j\mid a_j^{p_j}=1\rangle \to \langle a_i, a_j\mid a_i^{p_i}=a_j^{p_j}=1, \Delta_{a_i,a_j}^{r_{i,j}}=\Delta_{a_j,a_i}^{r_{i,j}}\rangle.
    \]
    If $\phi_{i,j}$ is injective, we define $\theta_{i,j}:=0$, otherwise $\theta_{i,j} := \pi/l_{i,j}$, where $l_{i,j}$ is half of the minimal syllable length of a non-trivial word in the kernel of $\phi_{i,j}$.
    
    From Lemma~\ref{lem:LengthReducedWordsRank2}, we know that a cyclically reduced non-empty word $w$ that is mapped to the identity has length $\geq 2r_{i,j}$. Therefore, we have that the natural maps $X\to G$, $Y\to G$ and $Z\to G$ are injective when $\frac{\pi}{{r_{1,2}}} + \frac{\pi}{{r_{2,3}}} + \frac{\pi}{{r_{1,3}}} \leq \pi$, or equivalently, $\frac{1}{{r_{1,2}}} + \frac{1}{{r_{2,3}}} + \frac{1}{{r_{1,3}}} \leq 1$. Hence, with this, we have proved that $\Sh_\Lambda \iso \Sh_\Gamma(\Lambda)$ whenever $\Lambda$ is a rank 2 induced subdiagram of $\Gamma$. 
    As all rank 2 Shephard groups are tame (a consequence of Lemma~\ref{lem:Rank2_properly_defined}), we get the same result for rank 1 induced subdiagrams.

    \smallskip \noindent \textbf{Case (d):}
    This case follows the construction used in case (c). Suppose now that there is at least one $r_{i,j}=\infty$. Without lost of generality, let us assume that $r_{1,2}=\infty$.
    Note that, when $r_{1,2}=\infty$, we have $\theta_{1,2}=0$. Hence $\phi_{1,2}$ is injective and $\theta_{1,2}+\theta_{2,3} + \theta_{1,3} = \theta_{2,3}+\theta_{1,3}\leq \pi$, for all $r_{i,j}\in \mathbb{N}_{\geq 2} \cup \{\infty\}$. Indeed, if both $\theta_{2,3}\neq \infty \neq \theta_{1,3}$, by Lemma~\ref{lem:LengthReducedWordsRank2}, we have that for all $r_{i,j}\in \mathbb{N}_{\geq 2}$, with $(i,j)\in\{(2,3),(1,3)\}$, $\frac{1}{r_{2,3}}+\frac{1}{r_{1,3}}\leq 1$. If $r_{2,3}$ and/or $r_{1,3}$ are $\infty$, we follow the same arguments and equally prove the injectivity of the natural maps $X\to G$, $Y\to G$ and $Z\to G$. As in case (c), we have $\Sh_\Gamma(\Lambda) = \Sh_\Lambda$ for any $\Lambda$ induced subdiagram of $\Gamma$.
%
%\textbf{Case (e):} In this case, we cannot use the triangle of groups, unless $k=\infty$, which is covered by case $(d)$. Hence, we will assume then $k\neq \infty$. We will suppose, without lost of generality, that $(r_{1,2},r_{2,3},r_{1,3})=(2,2,k)$ with $k\in \mathbb{N}_{\geq 2}$.
%Note that $r_{i,j}=2$ implies that $a_ia_j=a_ja_2$, i.e. in our case $a_2$ commutes with $a_1$ and $a_3$.
%Additionally, we have that $$\Sh_\Gamma \iso \langle a_2\mid a_2^{p_2}=1\rangle \times \langle a_1,a_3\mid a_1^{p_1}=a_3^{p_3}=1, \Delta_{a_1,a_3}^{r_{1,3}}=\Delta_{a_3,a_1}^{r_{1,3}}\rangle = \langle a_2\mid a_2^{p_2}=1\rangle \times \Sh_\Lambda,$$ where $\Lambda$ is the induced subgraph of $\Gamma$ with $V(\Lambda)=\{a_1,a_3\}$.
%By properties of direct products and the fact that $\Sh_\Lambda$ is tame, we have that $\Sh_\Gamma$ is also tame.
\end{proof}

A reducible rank $3$ Shephard group is tame by  Proposition~\ref{prop:tame_product} and Theorem~\ref{thm:Shephard_IP_rank2}.

\subsection{Finite Shephard Groups}

Following our characterization of the Buekenhout diagrams of rank $2$ finite Shephard groups, in this section we describe diagrams for higher ranks. In the following theorem, given an induced subgraph $\Lambda\subseteq \Gamma$, we we will denote by $\Lambda\setminus\{i\}$, for $i\in V(\Lambda)$, the induced subgraph obtained taking as vertex set $V(\Lambda)\setminus\{i\}$.

\begin{thm}\label{thm:FiniteShep_RC_FT}
    Let $\Gamma$ be a Shephard presentation graph.
    If $\Sh_\Gamma$ is finite, then $\Sigma(\Sh_\Gamma)$ is a residually-connected coset geometry on which $\Sh_{\Gamma}$ acts flag-transitively.
\end{thm}
\begin{proof}
    If $\Sh_\Gamma$ is a finite Shephard group of rank $2$, this is trivial.
    Additionally, this also true if $\Sh_\Gamma$ is a finite Coxeter group.
    Hence, consider that $\Sh_\Gamma$ is a finite Shephard group but not a Coxeter group.
    These groups are given in Table~\ref{table:finiteHigher}.
    
    From Theorem~\ref{thm:Shephard_finite_tame} and Theorem~\ref{thm:Tame_imply_IP}, we know that $\Sh_\Gamma$ satisfies the intersection property.
    Note that, by satisfying this property, we have that, for any induced subgraph $\Lambda\subseteq \Gamma$, $\Sh_\Gamma(\Lambda) = \langle \Sh_{\Gamma}(\Lambda\setminus\{i\})\mid i \in V(\Lambda)\rangle$.
    Moreover, since $\Gamma$ is a linear graph, there is a linear ordering of the generators of $\Sh_\Gamma$ such that for every non-consecutive pair of generators $\{a_i,a_j\}$ commute. Consequently, for every non-consecutive pair of generators $\{a_i,a_j\}$, the coset geometry $\Sigma(\langle a_i,a_j\rangle,(\langle a_i\rangle,\langle a_j\rangle)$ is going to be a generalized digon. 
    %One can then conclude that the diagram of $\Sigma(\Sh_\Gamma,(\Sh_\Gamma(\Gamma\setminus\{i\})_{i\in V(\Gamma)})$ will be a linear diagram.
    We can then use ~\cite[Theorem 2.6.4]{buekenhout2013diagram}, and conclude that $\Sigma(\Sh_\Gamma)$ is a residually-connected coset geometry and that $\Sh_\Gamma$ acts flag-transitively on it.
\end{proof}

\begin{table}[t]
\begin{tabular}{|c|c|c|}
\hline
\begin{tabular}[c]{@{}c@{}}Shephard\\ group\end{tabular} & Shephard presentation graph & Buekenhout diagram\\ \hline
$B_n(p,2)$ & $\xymatrix@-1pc{*{\bullet}\ar@{-}[rr]^4^(0.01)p && *{\bullet}\ar@{-}[rr] && *{\bullet}\ar@{-}[r] & {\cdots}\ar@{-}[r] &  *{\bullet}\ar@{-}[rr] && *{\bullet} }$ & 
$$\xymatrix@-1pc{*{\circ}\ar@{-}[rr]^{4,4,4}_(0.01){p-1} && *{\circ}\ar@{-}[rr]_(0.01){1}^{3,3,3} && *{\circ}\ar@{-}[r]_(0.01){1} & {\cdots}\ar@{-}[r] &  *{\circ}\ar@{-}[rr]_(0.01){1}_(0.99){1}^{3,3,3} && *{\circ} }$$ \\ \hline
$B_3(2,3)$ & 
$\xymatrix@-1pc{*{\bullet}\ar@{-}[rr]^4 && *{\bullet}\ar@{-}[rr]^(0.01)3^(0.99)3 && *{\bullet}}$
&
$\xymatrix@-1pc{*{\circ}\ar@{-}[rr]^{4,4,4}_(0.01){1} && *{\circ}\ar@{-}[rr]_(0.01){2}_(0.99){2}^{4,3,4} && *{\circ}}$
\\ \hline
$A_3(3)$ &
$\xymatrix@-1pc{*{\bullet}\ar@{-}[rr]^(0.01)3 && *{\bullet}\ar@{-}[rr]^(0.01)3^(0.99)3 && *{\bullet}}$
& 
$\xymatrix@-1pc{*{\circ}\ar@{-}[rr]_(0.01)2^{4,3,4} && *{\circ}\ar@{-}[rr]_(0.01)2_(0.99)2^{4,3,4} && *{\circ}}$
\\ \hline
$A_4(3)$  & 
$\xymatrix@-1pc{*{\bullet}\ar@{-}[rr]^(0.01)3 && *{\bullet}\ar@{-}[rr]^(0.01)3 && *{\bullet}\ar@{-}[rr]^(0.01)3^(0.99)3 && *{\bullet}}$
& 
$\xymatrix@-1pc{*{\circ}\ar@{-}[rr]_(0.01)2^{4,3,4} && *{\circ}\ar@{-}[rr]_(0.01)2^{4,3,4} && *{\circ}\ar@{-}[rr]_(0.01)2_(0.99)2^{4,3,4} && *{\circ}}$
\\ \hline
\end{tabular}
\caption{Buekenhout diagrams of the finite Shephard groups of rank $\geq 3$ that are not Coxeter groups. In the Shephard--Todd classification, they are $G(p,1,n)$, $G_{26}$, $G_{25}$, and $G_{32}$, respectively.}
\label{table:finiteHigher}
\end{table}

\begin{coro}
    The Buekenhout diagram of finite Shephard groups of rank $\geq 3$, which are non-Coxeter, are the ones given in Table~\ref{table:finiteHigher}.
\end{coro}

\begin{proof}
From Theorem~\ref{thm:FiniteShep_RC_FT}, we have that $\Sigma(\Sh_\Gamma)$ is a residually-connected coset geometry, and that $\Sh_\Gamma$ acts flag-transitively on it. This means that all rank $2$ residues of type $\{i,j\}$ are equal, for every distinct pair $i,j\in V(\Gamma)$.

Due to tameness of $\Gamma$ and the intersection property of $\Sh_\Gamma$, we have that each rank $2$ residue of type $\{i,j\}$ is exactly the coset geometry obtained for the finite rank $2$ Shephard group with presentation graph given by the induced subgraph of $\Gamma$ with vertex set $\{i,j\}$.
With this, one can then determine the Buekenhout diagram for the Shephard groups of Table~\ref{table:finiteHigher} from the Buekenhout diagrams given in Table~\ref{table:finiteR2}.   
\end{proof}

\section{Graph Product of Groups and Right-Angled Shephard Groups}

Armed with the main theorem of Section \ref{sec:higher}, stating that tame Shephard groups satisfy intersection property, we now investigate the well-behaved class of right-angled Shephard groups. Since right-angled Shephard groups are included inside Dyer groups, tameness, intersection property and flag-transitivity can be easily deduced for these groups, either from known results on Dyer groups, or results present in this paper.

Our goal here is to use a different approach. We will prove first that these three properties are present in any graph product of groups. As right-angled Shephard groups can be described as graph products, these results follow directly as a corollary.
First, let us define a graph product of groups.

\begin{definition}
    Given an unlabeled simple graph $\GG$ and a collection of groups $\C=\{G^v\mid v\in V(\GG)\}$, called the vertex groups, the \defn{graph product} $\GG \C$ is the group defined by
    \[
    (G^{v_1} \ast G^{v_2} \ast \cdots \ast G^{v_N}) / \langle gh=hg \mid g\in G^v, h\in G^k, \{v,k\} \notin E(\GG)\rangle,
    \]
    where $V(\GG) = \{v_1, v_2, \dotsc, v_N\}$.
\end{definition}

In other words, the group $\GG \C$ is the free product of all the vertex groups $G^v$, with additional relations imposing that elements of groups corresponding to non-adjacent vertices must commute. We alert the reader that this is not the usual definition of graph product of groups, where instead groups on adjacent vertices are considered to commute.
Indeed, in our definition,
if $\GG$ is the complete graph, then $\GG \C$ is the free-product of all $G^v$; while if $\GG$ has no edges, the group $\GG \C$ is the direct sum of the vertex groups.

Graph products can be used to define two important families of groups: right-angled Artin groups, where all vertex groups are infinite cyclic, and right-angled Coxeter groups, where all vertex groups are cyclic of order two. In this section, we will be interested in right-angled Shephard groups, which can be defined as graph products where all vertex groups are cyclic, without restriction on their sizes.

One important fact about graph products of groups is the existence of a unique reduced expression for an element $g\in \GG \C$, up to ``shuffling" (using the commuting relations)~\cite{Green1990,Genevois2019}. One direct result from this is very similar to the tameness of Shephard groups. Given an induced subgraph $\HH$ of $\GG$, the subgroup of $\GG \C$ generated by $\HH$, namely $\langle G^v \mid v\in V(\HH)\rangle$, is isomorphic to the graph product $\HH\LL$, where $\LL=\{G^v\mid v\in V(\HH)\}$.
Another result from Green is the fact that, for any $v\in \GG$, the subgroups of $\GG\C$ $\langle G^v\rangle$ and $\langle G^u\mid u\in V(\GG)\setminus\{v\}\rangle =\langle 1\rangle$ have trivial intersection (see~\cite[Lemma 3.3]{Green1990}).

\begin{thm}\label{thm:graphproduct}
    Let $\GG$ be an unlabeled simple graph, and let $\C=\{G^v\mid v\in V(\GG)\}$ be a collection of groups indexed by $\GG$.
    Then, the following hold.
    \begin{enumerate}
        \item For any pair of induced subgraphs $\HH,\KK\subseteq\GG$,we have
        $$\langle G^v\mid v\in V(\HH) \rangle \cap \langle G^v\mid v\in V(\KK)\rangle = \langle G^v\mid v\in V(\HH)\cap V(\KK) \rangle.$$
        \item The group $\GG \C$ acts flag-transitively on the coset incidence system $(\GG \C, (G_v)_{v\in V(\GG)})$, where $G_v = \langle G^u\mid u\in V(\GG)\setminus\{v\}\rangle$.
    \end{enumerate}
\end{thm}
\begin{proof}
    (a) One can derive this result in many different ways. Our proof will follow a similar approach as the one given in Theorem~\ref{thm:Tame_imply_IP}. Indeed, given an induced subgraph $\HH\subseteq \GG$, the fact that the subgroup of $\GG\C$ generated by $\HH$ is isomorphic to the graph product $\HH\LL$, with $\LL=\{G^v\mid b\in V(\HH)\}$ is equivalent to the definition of tameness. Hence, following the same construction given in Theorem~\ref{thm:Tame_imply_IP} we can derive the intersection property as defined in the statement of $(a)$.
    
    (b) We will prove this by induction on $n=|V(\GG)|$. For $n=1,2$ the theorem holds for rank one and two coset geometries since these are always flag-transitive.
    Suppose now, as the induction hypothesis, that for all graphs $\GG$ with less than $n$ vertices and all collections of groups $\C$, the graph product $\GG \C$ acts flag-transitively on the coset incidence system $(\GG\C, (G_v)_{v\in V(\GG)})$.
    Let $\GG$ be a graph with $n$ vertices and identify $V(\GG)$ with $\{0,1, \dotsc, n-1\}$. If $\GG$ has no edges on the $n$ vertices, the group $\GG \C$ is the direct product of the vertex groups and the theorem holds.
    We can thus assume that there exists a pair of distinct vertices $\{i,j\}$ that are connected by an edge of $\GG$.
    For any $J \subset \{ 0,1, \dotsc, n-1\}$, let $\GG_J$ be the subgraph of $\GG$ induced by $V \setminus J$ and $\C_J = \{(G_v)_{v \in V(\GG) \setminus J}\}$.
    By induction hypothesis, we know that $\GG_i \C_i$ and $\GG_j \C_j$ act flag transitively on their respective coset geometries.
    We now consider the free amalgamated product of $\GG_i \C_i *_{B}\GG_j \C_j$, where $B$ is the intersection of $\GG_i \C_i$ and $\GG_j \C_j$.
    By part (a), we know that this intersection is in fact $\GG_{\{i,j\}} \C_{\{i,j\}}$.
    This ensures that $\GG_i \C_i *_{B}\GG_j \C_j$ is actually isomorphic to $\GG \C$.
    We can then conclude that $\GG \C$ acts transitively on $(\GG \C, (G_v)_{v\in V(\GG)})$ by the results on coset geometries of amalgamated products in \cite{piedade2026group}.
\end{proof} 

We now come back to Shephard groups.
Recall a Shephard group is right-angled if all edges of the presentation graph are labeled by $\infty$ (Definition~\ref{def:RAShepG}).

In other words, two generators of a right-angled Shephard group either commute or generate a free group of rank two. In particular, it follows that a right-angled Shephard group is a graph product where all vertex groups are cyclic.

\begin{coro}
    Right-angled Shephard group are tame, satisfy the intersection property and act flag-transitively on their natural coset geometry.
\end{coro}

\begin{proof}
    Part (a) of Theorem \ref{thm:graphproduct} is a rephrasing of tameness to the context of graph products. We thus conclude that right-angled Shephard groups are tame. Hence, by Theorem \ref{thm:Tame_imply_IP}, it satisfies intersection property. Lastly, part (b) of Theorem \ref{thm:graphproduct} states that the group acts flag-transitively on its standard coset geometry.
\end{proof}

\section{Representations of Shephard Groups}

We start with the Burau representation associated with Artin--Tits groups (i.e., Shephard groups whose vertices are all labeled by $\infty$) following the description given by Bapat--Queffelec~\cite{BQ24}.
To do so, we have the $t$-deformed pairing
\[
\langle \alpha_i, \alpha_j \rangle = \begin{cases}
1 + t^2 & \text{if } i = j, \\
-2t \cos(\pi / m_{ij}) & \text{otherwise},
\end{cases}
\]
where by the usual conventions $\cos(\pi / \infty) = \cos(0) = 1$.
(We have performed the substitution $q \mapsto -t$ from~\cite{BQ24} to more closely match the Coxeter geometric representation; see~\cite[Sec.~4.2]{BB05}.)
This pairing is $t$-sesquilinear in the sense
\[
\langle tx, y \rangle = t^{-1} \langle x, y \rangle = \langle x, t^{-1} y \rangle.
\]

\begin{definition}
    Let $R$ be a commutative ring containing $t$ and $\cos(\pi / m_{ij})$ for all $i, j \in V(\Gamma)$.
    The \defn{Burau representation} is the action on $V_{\Gamma} := \bigoplus_{i \in V(\Gamma)} R \alpha_i$ given by
    \begin{equation}
        \label{eq:geometric_repr}
        a_i(\alpha_j) = \alpha_j - \langle \alpha_i, \alpha_j \rangle \alpha_i.
    \end{equation}
    We denote this by $\rho_{\Gamma} \colon \Sh_{\Gamma} \to GL(V_{\Gamma})$.
\end{definition}

\begin{remark}
    This is the case of $s = t$ in Squier~\cite{Squier88}, which we believe to be the first place these representations were introduced for general Artin groups.
\end{remark}

Next, we recall the \defn{geometric representation} of the  Shephard group $\Sh_{\Gamma}$ when $p_i < \infty$ for all $i \in V(\Gamma)$ following~\cite[Def.~2.1]{Goldman2024}, which is adapted from~\cite{Mostow80}.
In what follows, we will assume $R$ is large enough so the result of every computation is in $R$.
We begin by defining a symmetric form $(\ ,\ ) \colon V_{\Gamma} \times V_{\Gamma} \to R$ by the following three cases:
\begin{enumerate}
\item if $i \neq j$ and $m_{ij} \neq 2$, then
$
\displaystyle
(\alpha_i, \alpha_j) := -\sqrt{\frac{\cos(\pi/p_i - \pi/p_j) + \cos(2\pi/m_{ij})}{2\sin(\pi/p_i)\sin(\pi/p_j)}};
$
\item if $i \neq j$ and $m_{ij} = 2$, then $(\alpha_i, \alpha_j) := 0$;
\item $(\alpha_i, \alpha_i) := 1$.
\end{enumerate}
Next, we define the bilinear form $\langle\ ,\ \rangle \colon V_{\Gamma} \times V_{\Gamma} \to R$ by
\[
\langle \alpha_i, \alpha_j \rangle := (1 - \zeta_{p_i}) (\alpha_i, \alpha_j),
\]
where $\zeta_m := e^{2\pi \ii / m}$ is a primitive $m$th root of unity (with $\ii = \sqrt{-1}$ as usual).
Using this bilinear form, we define the geometric representation $\rho_{\Gamma} \colon \Sh_{\Gamma} \to GL(V_{\Gamma})$ by Equation~\eqref{eq:geometric_repr} as for the Burau representation.

The goal of this section is effectively to combine these two representations together.
Towards this goal, let us first consider the case when $p_i = p_j = p$ (which is necessarily true when $m_{ij}$ is odd) to show how we can derive the Burau representation from the geometric representation as a sort of $p \to \infty$ limit.
In this case, we only need to consider when $i \neq j$ and $m_{ij} \neq 2$, and hence
\begin{equation}
\label{eq:equal_param_pair}
( \alpha_i, \alpha_j ) = -\sqrt{\frac{1 + \cos(2\pi/m_{ij})}{2\sin^2(\pi/p)}}
= - \frac{\cos(\pi/m_{ij})}{\sin(\pi/p)},
\end{equation}
by the double angle formula and $\sin(\pi/x), \cos(\pi/x) \geq 0$ for all $x \geq 2$.
While taking the limit
\[
\lim_{p \to \infty} \frac{1 - e^{2\pi\ii/p}}{\sin(\pi/p)} = -2 \ii
\]
means we do get a valid representation, where the $(i,j)$ submatrix of this representation being
\[
a_i \mapsto \begin{bmatrix} 1 & 2 \ii \cos(\pi/m_{ij}) \\ 0 & 1 \end{bmatrix},
\qquad\qquad
a_j \mapsto \begin{bmatrix} 1 & 0 \\ 2 \ii \cos(\pi/m_{ij}) & 1 \end{bmatrix}.
\]
Yet, it becomes the trivial representation for the rank $1$ case since $e^{2\pi \ii/p} - 1 \to 0$ as $p \to \infty$ (as indicated by the diagonals of the above matrices).
Instead, let us consider $\ii t = e^{\pi\ii/p} = \frac{d}{dP} e^{\ii t P} \rvert_{P=0}$, for $P = t^{-1} \pi p^{-1}$, as an infinitesimal generator of the circle group $\RR / \ZZ$, in which case we replace $\zeta_{p_i} = -t^2$ and compute
\[
\frac{1 - e^{2\pi\ii/p}}{\sin(\pi/p)} = \frac{1 + t^2}{\frac{t + t^{-1}}{2}} = 2 t,
\]
which means we recovered the Burau representation.

Therefore, to get a representation $\rho_{\Gamma} \colon \Sh_{\Gamma} \to GL(V_{\Gamma})$ of the general Shephard group $\Sh_{\Gamma}$, we only need to consider the case when, without loss of generality, $p_j < p_i = \infty$.
While there are no issues when taking $\lim_{p_i\to\infty} \cos(\pi/p_i - \pi/p_j) = \cos(\pi/p_j)$, we still need to use our infinitesmial generators instead of directly evaluating this limit.
Hence, applying the angle subtraction formula and taking the substitutions
\[
\sin(\pi/p_i) \mapsto \frac{t + t^{-1}}{2},
\qquad\qquad
\cos(\pi/p_i) \mapsto \frac{t - t^{-1}}{2} \ii,
\]
we obtain
\begin{align*}
(\alpha_i, \alpha_j) % & = -\sqrt{ \frac{\cos(\pi/p_j)(t^2 - t^{-1})\ii + \sin(\pi/p_j)(t + t^{-1}) + 2\cos(2\pi/m_{ij})}{2 (t + t^{-1})\sin(\pi/p_j)}}
% \\ & = -\sqrt{ \frac{t e^{-\ii \pi / p_j} \ii + (t e^{-\ii \pi / p_j} \ii)^{-1} + 2\cos(2\pi/m_{ij})}{2 (t + t^{-1})\sin(\pi/p_j)}}
 & = -\sqrt{ \frac{t \zeta_{2p_j}^{-1} \ii + (t \zeta_{2p_j}^{-1} \ii)^{-1} + 2\cos(2\pi/m_{ij})}{2 (t + t^{-1})\sin(\pi/p_j)}}.
\end{align*}
(Since $\ii = \zeta_4$, we could simplify $\zeta_{2\pi_j}^{-1} \ii$ into a power of a single root of unity.)
Equivalently, this comes from taking the substitution
\begin{equation}
\label{eq:inf_limit}
p_i^{-1} = -\ii \pi^{-1} \log \ii t,
\end{equation}
which is how we obtain the infinitesimal generator.
Hence, our $(i,j)$ submatrix of this representation becomes
\[
a_i \mapsto \begin{bmatrix} -t^2 & (1 + t^2) (\alpha_i, \alpha_j) \\ 0 & 1 \end{bmatrix},
\qquad\qquad
a_j \mapsto \begin{bmatrix} 1 & 0 \\ (1 - \zeta_{p_j}) (\alpha_j, \alpha_i) & \zeta_{p_j} \end{bmatrix}.
\]

\begin{thm}
    The map $\rho_{\Gamma} \colon \Sh_{\Gamma} \to GL(V_{\Gamma})$ is a representation of the Shephard group $\Sh_\Gamma$.
\end{thm}

\begin{proof}
Since the Coxeter relations $\Delta_{a_i,a_j}^{r_{ij}}$ hold independently of $p_i$ and $p_j$, we can simply apply the substitution in~\eqref{eq:inf_limit} to the matrices of the relation.
\end{proof}

Likewise, we call $\rho_{\Gamma}$ the \defn{geometric representation} of the Shephard group $\Sh_{\Gamma}$ (which includes the Burau representation for Artin--Tits groups).

Now we focus on the rank $2$ cases.

\begin{lemma}
    \label{lem:simple_geom_repr}
    For an irreducible rank $2$ Shephard group, the geometric representation is irreducible.
\end{lemma}

\begin{proof}
    We prove the case when $p_i, p_j < \infty$ as the cases for when one or both are infinite is similar.
    Since the geometric representation is two dimensional, we only need to show there is not an invariant one dimensional subspace.
    As $\rho_{\Gamma}(a_1)$ is an upper triangular matrix, its eigenvectors are $(1,0)$ and $(1, -\langle \alpha_i, \alpha_j \rangle)$.
    Similarly, $\rho_{\Gamma}(a_2)$ being lower triangular matrix means its eigenvectors are $(0,1)$ and $(-\langle \alpha_j, \alpha_i \rangle, 1)$.
    Therefore, there is no simultaneous eigenspace for these matrices, and the geometric representation is irreducible.
\end{proof}

Our proof is effectively shows the bilinear form $\langle\ ,\ \rangle$ is non-degenerate as its Gram matrix for the basis $\{\alpha_i, \alpha_j\}$ is
\[
\begin{bmatrix}
1 & (1 - \zeta_{p_i}) (\alpha_i, \alpha_j) \\
(1 - \zeta_{p_j}) (\alpha_i, \alpha_j) & 1
\end{bmatrix}
\quad \xrightarrow{\quad\det\quad} \quad
1 - (1 - \zeta_{p_i}) (1 - \zeta_{p_j}) (\alpha_i, \alpha_j)^2 \neq 0,
\]
where we could then invoke~\cite[Lemma~2.3]{Goldman2024}.

\begin{thm}
    \label{thm:rank2_faithful}
    Consider the rank $2$ Shephard group $\Sh_\Gamma = \Sh(p_i,r,p_j)$ with $p_i, p_j < \infty$, then the geometric representation $\rho_{\Gamma}$ is faithful if and only if
    \[
    \displaystyle \frac{1}{p_i} + \frac{1}{p_j} + \frac{2}{r} > 1.
    \]
\end{thm}

\begin{proof}
From~\cite{Mostow80}, the group is finite, the geometric representation defined above is faithful, and the group is finite if and only if $\frac{1}{p_i} + \frac{1}{p_j} + \frac{2}{r} > 1$ by~\cite[Thm.~2.7]{Goldman2024}.

Now suppose $p_i, p_j < \infty$ but $\frac{1}{p_i} + \frac{1}{p_j} + \frac{2}{r} \leq 1$.
First, we consider the case $r = 2m$ is even, in which case by~\cite[Lemma~3.2]{2ShepKatherine}, the center of $\Sh_{\Gamma}$ is generated by $c = (a_i a_j)^m = (a_j a_i)^m$ and is infinite cyclic.
Since $c$ is in the center and that the geometric representation is irreducible by Lemma~\ref{lem:simple_geom_repr}, Schur's Lemma implies that $\rho_{\Gamma}(c)$ must be a scalar matrix $\lambda I$ (by considering it as a $\Sh_{\Gamma}$-linear map from $V_{\Gamma} \to V_{\Gamma}$).
Next, we compute the matrix product
\begin{align*}
\rho_{\Gamma}(a_i a_j) & = \begin{bmatrix}
\zeta_{p_i} + (1 - \zeta_{p_i}) (1 - \zeta_{p_j}) (\alpha_i, \alpha_j)^2 & (1 - \zeta_{p_i}) \zeta_{p_j} (\alpha_i, \alpha_j) \\
(1 - \zeta_{p_j}) (\alpha_i, \alpha_j) & \zeta_{p_j}
\end{bmatrix},
%\\
%\rho_{\Gamma}(a_j a_i) & = \begin{bmatrix}
%\zeta_{p_i} & (1 - \zeta_{p_i}) (\alpha_i, \alpha_j) \\
%\zeta_{p_i} (1 - \zeta_{p_j}) (\alpha_i, \alpha_j) & \zeta_{p_j} + (1 - \zeta_{p_i}) (1 - \zeta_{p_j}) (\alpha_i, \alpha_j)^2
%\end{bmatrix},
\end{align*}
and by adding $(\zeta_{p_j} - 1) (\alpha_i, \alpha_j)$ times the second row to the first, we obtain
\[
E_{12}\bigl((\zeta_{p_j}-1) (\alpha_i, \alpha_j)\bigr) \rho_{\Gamma}(a_i a_j) = \begin{bmatrix}
\zeta_{p_i} & 0 \\
(1 - \zeta_{p_j}) (\alpha_i, \alpha_j) & \zeta_{p_j}.
\end{bmatrix}
\]
Therefore, we see the determinant of $\rho_{\Gamma}(a_i a_j)$ is $\zeta_{p_i} \zeta_{p_j}$, and hence the determinant of $\rho_{\Gamma}(c)$ is $(\zeta_{p_i} \zeta_{p_j})^m$, which implies that $\rho_{\Gamma}(c) = (\zeta_{p_i} \zeta_{p_j})^{m/2} I$.
This is a matrix of finite multiplicative order, and the geometric representation cannot be faithful in this case.

Now let us consider the case $p_i = p_j = p < \infty$, and we note that $\Sh(p,r,p)$ is isomorphic to an index $2$ subgroup of $\Sh(p,2r,2) = \langle b_i, b_j \rangle$ generated by $a_i = b_i$ and $a_j = b_j b_i b_j$~\cite[Lemma~3.3]{2ShepKatherine}.
Therefore, we have $(a_i a_j)^k = (b_i b_j b_i b_j)^k = (b_i b_j)^{2k}$, and from~\cite[Lemma~3.2]{2ShepKatherine}, we have $(b_i b_j)^{2r}$ is in the center.
From the previous argument, we have that $\rho(a_i a_j)$ has finite multiplicative order, but the subgroup generated by $(b_i b_j)^{2k}$ is infinite.
Hence, the representation is not faithful.
\end{proof}

From~\cite[Thm.~2]{Squier88}, when $p_i = p_j = \infty$, then the geometric representation (which is the Burau representation) is faithful in the rank $2$ case.
Thus, we have the following outstanding case, for which we believe faithfulness still holds.

\begin{conjecture}
    For $p_i < p_j = \infty$, then the geometric representation $\rho_{\Gamma}$ is faithful.
\end{conjecture}

We note that we can follow the argument in Theorem~\ref{thm:rank2_faithful}, but the resulting matrix has infinite order as the scalar contains a power of $t$ (assuming that $t$ is not a root of unity).

\subsection*{Acknowledgments}

The authors thank Katherine Goldman for many useful discussions regarding tameness of Shephard groups and their representation; in particular, we would like to thank for the idea of proving Theorem~\ref{thm:rank2_faithful} by looking at elements in the center of $\Sh(p,2r,q)$.

CAP is funded by an Action de Recherche Concert\'{e}e -- ARC -- from the Communaut\'{e}e Fran\c{c}aise Wallonie-Bruxelles. 
TS was partially supported by Grant-in-Aid for JSPS Fellows 21F51028 and JSPS KAKENHI Grant Number JP23K12983.
PT is a Postdoctoral Researcher of the Fonds de la Recherche Scientifique -- FNRS.

\bibliographystyle{siam}
\bibliography{refs}

\end{document}